\RequirePackage[l2tabu,orthodox]{nag} 

\documentclass[a4paper,12pt,english]{article}

\usepackage{setspace}

\usepackage{authblk}

\usepackage{ulem}

	\usepackage{lmodern}				

	\usepackage[T1]{fontenc}			

	\usepackage[utf8x]{inputenc}		

\usepackage{graphicx} 
\usepackage{tikz}

	\usepackage{babel}					

\usepackage{amsmath}		
\usepackage{amsthm}	
\usepackage{amssymb}		

\RequirePackage{bbm}					

\RequirePackage{stmaryrd}

	\usepackage{pstricks}
	\usepackage{hyperref}
	\usepackage{cleveref}

\usepackage{geometry}

\usepackage{algorithm}
\usepackage{algorithmic}

\usepackage{subcaption}

\usepackage{url}

\newcommand{\ensnombre}[1]{\mathbb{#1}}%

\newcommand{\N}{\ensnombre{N}}

\newcommand{\R}{\ensnombre{R}}

\newcommand{\abs}[1]{\left \lvert #1 \right \rvert}

\newcommand{\norme}[1]{\left \lVert #1 \right \rVert}

\newcommand{\defeq}{\mathrel{\mathop:}=}

\newcommand{\eqdef}{\mathrel{=}:}

\newcommand{\ind}{\mathbbm{1}}

\renewcommand{\P}{\mathbb{P}}

\newcommand{\E}{\mathbb{E}}

\def \Var{\hbox{{\textrm{Var}}}}

\def \Cov{\hbox{{\textrm{Cov}}}}

\definecolor{vert_fonce}{rgb}{0,0.75,0.25}
\definecolor{vert_pire}{rgb}{0,0.70,0.30}
\definecolor{vert_plus}{rgb}{0,0.60,0.40}
\definecolor{rouge_plus}{rgb}{0.70,0.20,0.10}
\definecolor{bleu_plus}{rgb}{0.30,0,0.70}
\definecolor{green_moi}{RGB}{0,100,0}
\definecolor{forestgreen}{rgb}{0,0.6,0}

\theoremstyle{plain}

\newtheorem{thm}{Theorem}[section]

\newtheorem{cor}[thm]{Corollary}

\newtheorem{lem}[thm]{Lemma}

\newtheorem{prop}[thm]{Proposition}

\theoremstyle{remark}

\newtheorem{rem}[thm]{Remark}

\usepackage{setspace}
\usepackage{comment}

\makeatletter
\def\input@path{{../Soumission/}}
\makeatother

\graphicspath{{../Soumission/}}

\renewcommand{\geq}{\geqslant}
\renewcommand{\leq}{\leqslant}

\usepackage{authblk}
\title{
Sensitivity analysis from a single input/output sample}
\author[1]{S\'ebastien Da Veiga}
\author[2]{Fabrice Gamboa}
\author[3]{Thierry Klein}
\author[4]{Agn\`es Lagnoux}
\author[5]{Cl\'ementine Prieur}
\affil[1]{\footnotesize{Universit\'e de Rennes, Ensai, CNRS, CREST - UMR 9194, F-35000 Rennes, France}}
\affil[2]{\footnotesize{Institut de Math\'ematiques de Toulouse and ANITI; UMR5219. Universit\'e de Toulouse; CNRS. UT3, F-31062 Toulouse, France and ANITI Toulouse France}}
\affil[3]{\footnotesize{Institut de Math\'ematiques de Toulouse; UMR5219. Universit\'e de Toulouse; ENAC - Ecole Nationale de l'Aviation Civile , Universit\'e de Toulouse, France}}
\affil[4]{\footnotesize{Institut de Math\'ematiques de Toulouse and ANITI; UMR5219. Universit\'e de Toulouse; CNRS. UT2J,  F-31058 Toulouse, France and ANITI Toulouse France}}
\affil[5]{\footnotesize{Universit\'e de Grenoble Alpes, CNRS, Inria, Grenoble INP, LJK, Grenoble, France}}
\begin{document}
\maketitle


\begin{abstract}
	The main objective of this  paper is to estimate optimally Sobol' indices at any order when a unique input/output i.i.d.\ sample is available. 
	Our approach stands on  three main ingredients: semi-parametric estimation theory, high-order kernel estimation 
	(inspired by the paper \cite{doksum1995nonparametric}), and  mirror-type transformations as introduced in \cite{bertin2020adaptive,pujol2022nonparametric}. We propose two different estimators. We prove that these estimators are asymptotically normal and efficient. Furthermore, we illustrate their numerical  properties on standard examples. 
\end{abstract}



\section{Introduction} \label{sec:intro}
In this paper, {we consider a random variable $Y$ depending on a random vector $V=(V_1,\ldots,V_p)$ through the relation $Y=G(V_1,\ldots, V_p)$ where the function $G$ is defined on a compact subset of $\R^p$, $p \geq 2$ and is real-valued. Then we}
 focus on the estimation of the ratio
\begin{equation*}
S^X=\frac{\Var(\E[Y|X])}{\Var(Y)}=\frac{\E[\E[Y|X]^2]-\E[Y]^2}{\Var(Y)}
\end{equation*}
where $X$ is a $d$-dimensional {vector formed by $d$ coordinates of} $V$ {($d$ being a non zero integer smaller than $p$)}. The interest in this ratio dates back to Karl Pearson as part of analysis of variance (see, e.g., \cite{kenstu67}), where it coincides to the square of the correlation ratio (measure of the relationship between the statistical dispersion within individual categories and the dispersion across a whole population or a sample). It was later revisited by Sobol' in the framework of sensitivity analysis \cite{sobol1993,sobol2001global}, under the name of \textit{closed Sobol' index}. Interestingly, it also appeared in the statistics literature under different names: \textit{measure of explanatory power of covariates} or \textit{nonparametric coefficient of determination} by \cite{doksum1995nonparametric}, or the \textit{residual variance} (focusing only on the numerator) which provides a lower bound for the performance of any regression function estimator in multivariate problems, as noticed in \cite{devroye2018nearest}. Finally, recent studies have emphasized the connection between machine learning feature importance measures and this ratio (see \cite{borgonovo2024total} and references therein). In particular, it is the building block of the \textit{Shapley values} which are now highly popular importance measures in machine learning for explainability \cite{owen2014sobol}. {It can also be used for recursive feature elimination \cite{benard2022mean}.}

\medskip

Here, without loss of generality we will adopt the notation and terminology from the field of sensitivity analysis, where this ratio is called the closed Sobol' index, hence the notation $S^X$ above. This research area emerged in the computer experiments community, which focuses on complex computer models to simulate and analyze natural systems in physics,
engineering and other fields. These models usually depend on many
input variables, and it is thus crucial to understand which input parameter or which set
of input parameters have an influence on the output. This is exactly the purpose of sensitivity analysis, which has become an essential tool for systems modelling and policy support (see, e.g., \cite{razavi:hal-03518573}). Global sensitivity analysis methods consider the input vector as random and propose a measure of the influence, in terms of output fluctuations, of each subset of its components. We refer to the seminal book \cite{saltelli-sensitivity} for an overview on global sensitivity analysis or to \cite{da2021basics} for a synthesis of recent trends in this field. 
Among the different measures of global sensitivity analysis, variance-based measures are probably the most commonly used. More precisely, for the output $Y$ of a computer code $G$, one of the most common measure of the sensitivity of $Y$ with respect to a {vector $X$} is precisely the closed Sobol' index $S^X$ {as defined at the very beginning of this section}.

\medskip

In recent years a myriad of different estimators have been proposed, see \cite[Chapter 4]{da2021basics} for a full review. A practical statistical method widely used to construct estimators is based on spectral analysis of the input/output functional relationship and Parseval formulae. We refer to \cite{sudret2008global} for a basic description of the method and to \cite[Chapter 4]{da2021basics} for more recent references. It should be noted that the asymptotic properties of these methods  have been little studied, as they are based on the theory of non-linear (quadratic) functional estimation.  
Away from the spectral methods, two families of estimation methods are of particular interest because of the possibility of studying the asymptotic properties of the estimators produced.

\medskip

The first family of methods is based on the so-called Pick-Freeze (PF) design of experiments. The basic idea is to evaluate $G$ repeatedly for input values, fixing those of the vector $X$, and then to calculate a Monte Carlo estimate of $S^X$ from this particular sampling. The main advantage is that only minimal assumptions are required to derive consistency and central limit theorems under the assumption of independent covariates. In particular, assumptions of integrability but not regularity on $G$ are necessary (see for example \cite{pickfreeze, janon2012asymptotic}). As an illustration, to estimate at rate $\sqrt{n}$ a single Sobol' index, one needs a design of experiments of size $2n$.  This implies that estimating all the $p$ first-order indices corresponding to $d=1$ involves a sample of size $(p+1)n$; a $(p+2)n$ sample is required to estimate the full set of first-order, second-order, and total indices. This cost can be reduced to $2n$ \cite{tissot2015randomized,gilquin2019making} if the aim is only to estimate first-order and second-order indices. Recently, the authors in \cite{borgonovo2024total} proposed to apply a weighting factor (called density quotient) to PF estimators to handle dependent covariates in settings in which
both probability distributions of $V$ and $X$ are known. However, the main drawback of the PF estimation procedure is that it requires a specific experimental design. In particular, it cannot be used in the case where we only have a classical i.i.d.\  sample of $n$ input/output observations. 
 
 \medskip

The second class of methods relies on local averaging and tackles this limitation. Among them are kernel estimators, which have been thoroughly studied for the case $d=1$ \cite{da2009local,da2008efficient,plischke2020fighting,solis2021non,heredia2021nonparametric} with central limit theorems and asymptotic efficiency as soon as $G$ satisfies regularity assumptions. Closely related are nearest neighbor approaches, which have been studied by several authors (see, e.g., \cite{devroye2003estimation,liitiainen2008nonparametric,
liitiainen2010residual,devroye2013strong,
gyorfi2015asymptotic,devroye2018nearest}). For instance in \cite{devroye2018nearest}, the authors propose a plug-in estimator with statistical consistency for any $d$ and a central limit theorem with rate $\sqrt{n}$ for $d\leq 3$ if, again, regularity assumptions hold. In parallel, \cite{broto2020variance} consider a variant which is consistent for any $d$ but no rate of convergence is provided. When $d=1$, a central limit theorem for estimators based on ranks (i.e.\ nearest neighbors on the right) is also proved in \cite{GGKL20}.

\medskip

In a nutshell, on the one hand, the class of PF estimators allows to estimate Sobol' indices at rate $\sqrt{n}$ for any $d$ with minimal assumptions on $G$ but requires a sample with highly specific structure. On the other hand, local-averaging estimators are base on {a single input-output $n$-sample} {(standard Monte-Carlo)} but need some regularity assumptions on the model and  $\sqrt{n}$-parametric rate of convergence is proved only for $d\leq 3$ (nearest neighbors method). In this paper, we propose and study a more general local-averaging estimator and show that the rate $\sqrt{n}$ holds for any $d$.

\medskip

In our work, we build an optimal estimator. In that view, we correct a plug-in estimator using a high-order kernel and a mirror-type transformation. The correction improves the initial estimator by correcting the bias (see, e.g., \cite{pfanzagl1982lecture,
newey1998undersmoothing}), in line with techniques used to improve an estimator in the frame of semi-parametric inference. We refer to \cite[Chapter 25]{van2000asymptotic} for an exhaustive overview of this theory. As a matter of fact, in our frame, the nonparametric part is the unknown regression function.  It is handled by using a very specific kernel estimator. More precisely, on the one hand,
we make use of high-order kernels in order to obtain $\sqrt{n}$ rate of convergence, following recent approaches of nonparametric regression (see for example \cite{tsybakov2009nonparametric}). Note that a high-order kernel has a non-zero negative part thus is not a probability density. Its use allows to {cancel the first terms in the bias expansion (see the details of such computations in \eqref{eq:bias_kill}}. On the other hand, to remove boundary effects inherent to the kernel estimation procedure, we adapt recent mirror-type transformations \cite{bertin2020adaptive,pujol2022nonparametric}. We propose and study two variants for the mirror transformation. In both cases, we show that the proposed estimators satisfy a central limit theorem with  the optimal rate and minimal variance. Up to our knowledge, this general optimality result is new in the frame of Sobol' index estimation based on {a single input-output $n$-sample}. Notice that the order one case is tackled in \cite{da2008efficient} by using another estimator based on classical kernel estimation jointly with a spectral approach. Note finally that our results are closely related to the ones developed in \cite{doksum1995nonparametric} where the boundary issue is tackled using truncation leading to the estimation of a pruned version of $S^X$.

\medskip

The paper is organized as follows. The setting and the notation are introduced in Section \ref{sec:setting}. Then, considering two different mirror-type transformations, we present in Section \ref{sec:estimation_procedure} two kernel-based regression estimators. Section \ref{sec:main_compact} is devoted to the statement of our main results, namely central limit theorems and asymptotic efficiency. A comparison study of asymptotic estimation variance with existing methods is provided in Section \ref{sec:compar}. Section \ref{sec:num} presents several numerical illustrations of our estimation procedure. The proofs are postponed to Appendix \ref{app:proof}. Moreover, extended numerical studies discussing numerical stability of high-order kernels can be found as an additional file on the following link \url{https://hal.science/hal-04052837}.

\section{Setting}\label{sec:setting}

\subsection{Model and purpose}\label{sec:model}

The output $Y$ is obtained from  a regression model (generally computed by a numerical code or a machine learning estimated model), and is regarded as a function $G$ of random inputs $(V_i)_{i=1,\dots,p}$ such that $Y=G(V_1,\ldots, V_p)$ where the function $G$ is defined on a compact subset of $\R^p$, $p \geq 2$ and is real-valued. {The random vector formed by the $p$ random input variables is denoted by $V$.}
Global sensitivity analysis is performed using Hoeffding decomposition \cite{Hoeffding48} for $Y\in \mathbb L^2(\R)$ which leads to standard Sobol' indices introduced below \cite{sobol1993} assuming that the $V_i$'s are independent random variables. However, such independence assumption is only required for practical interpretation of the decomposition, but is not at all necessary for estimation purposes and asymptotic guarantees. Thus, in the rest of the paper, we do not work with such strong assumption.

 \medskip
 
For any $u \subsetneq \{1,\dots,p\}$, $u\neq \emptyset$ we denote $X=(V_i)_{i\in u}$ a group of inputs with cardinality $\abs{u}=d <p$. The \textit{closed Sobol' index} of output $Y$ associated to the $d$-dimensional vector $X$ is defined as
 \begin{equation}\label{def:sobol}
S^X=\frac{\Var(\E[Y|X])}{\Var(Y)}=\frac{\E[\E[Y|X]^2]-\E[Y]^2}{\Var(Y)}.
\end{equation}
For the specific case $u=\{i_0\}$ with $i_0\in \{1,\dots,p\}$,  that is for $X=V_{i_0}$,   $S^{i_0}=S^X$ is the first-order Sobol' index associated to the input $V_{i_0}$. In addition, the \textit{total Sobol' index} associated to $V_{i_0}$ writes 
{\begin{equation}\label{def:sobol_total}
    S^{\text{tot},i_0}=1-S^{V_{\sim i_0}}=1-S^{\sim i_0}
\end{equation} }
where  $V_{\sim i_0}= (V_1,\cdots, V_{i_0-1},V_{i_0+1}, \cdots, V_p)$.

\medskip

Our goal is to estimate $S^X$ from a $n$-sample $(X_j,Y_j)_{j=1,\dots,n}$ of the pair $(X,Y)$ with joint distribution $\mathcal{P}$. Since $\E[Y]$ and $\Var(Y)$ can naturally be estimated with classical empirical moments, we focus here on the estimation in \eqref{def:sobol} of $T=\E[\E[Y|X]^2]$.

\subsection{Notation and assumptions}\label{sec:nota}

In this section, we give the general assumptions made in the paper.

\renewcommand{\theenumi}{\textbf{($\mathcal A$\arabic{enumi})}}
\renewcommand{\labelenumi}{\theenumi}
\begin{enumerate}
\setcounter{enumi}{0}
\item\label{hyp:domain}  The support  of $(V_1,\ldots,V_p)$ is $[0,1]^p$. { We assume that $X$ is absolutely continuous with respect to the Lebesgue measure on $[0,1]^d$ with density function $f_X$. }
\end{enumerate}
\renewcommand{\theenumi}{\alph{enumi}}
\renewcommand{\labelenumi}{\theenumi.}

For simpler notation here we only consider $[0,1]^p$ but results can be readily extended for support of the form $[B_1,C_1]\times \cdots \times  [B_p,C_p]$ where $B_i<C_i$ for all $i \in \{1,\dots,p\}$.

\renewcommand{\theenumi}{\textbf{($\mathcal A$\arabic{enumi})}}
\renewcommand{\labelenumi}{\theenumi}
\begin{enumerate}
\setcounter{enumi}{1}
\item\label{hyp:densityb}
$\exists \, c >0$ such that $\inf_{x\in [0,1]^d} f_X(x)\geqslant c$. 

\item\label{hyp:bounds} $\E[Y^4]<\infty$ and $\sigma^2(x)=\Var(Y\vert X=x)$ is bounded on $[0,1]^d$. 
\end{enumerate}
\renewcommand{\theenumi}{\alph{enumi}}
\renewcommand{\labelenumi}{\theenumi.}

Note that under Assumption \ref{hyp:bounds} the regression function $m(\cdot)=\E[Y|X=\cdot]$ exists and is bounded on the support of $X$ by Cauchy-Schwartz inequality. Furthermore, we define the function $g\defeq f_Xm$.

\medskip

For any integer $d$, we introduce the following multi-index notation. For any $\beta =(\beta_1, \ldots , \beta_d) \in (\mathbb{R}_{>0})^d$,  let  $\lfloor \beta\rfloor$ be the largest integer stricly lower than $\beta$: $\lfloor \beta \rfloor = (\lfloor \beta_1 \rfloor, \dots, \lfloor \beta_d \rfloor) \eqdef \gamma \in \N^d$.
In addition, we introduce 
\[
 |\gamma| = \gamma_1 +\dots + \gamma_d \, , \;  \gamma !=\gamma_1 !\dots  \gamma_d!, \textup{ and }  
 v^\beta = v_1^{\beta_1}\dots v_d^{\beta_d} \textup{ for any }v\in \R^d \, .
\]


For $\alpha >0$ let $\mathcal{C}^{\alpha}([0,1]^d)$ be the set of functions $\phi \colon [0,1]^d \to \R$ that have derivatives up to order $\lfloor \alpha \rfloor$ and for which partial  derivative of order $\lfloor \alpha \rfloor$ is $\alpha - \lfloor \alpha \rfloor$-H\"older. Namely, there exists $C_{\phi}>0$ such that, for any $x$ and $x'\in [0,1]^d$, one has 
\begin{align}\label{eq:holder}
\Bigl\vert \frac{\partial^\beta \phi}{\partial x^\beta}(x) - \frac{\partial^\beta \phi}{\partial x^\beta}(x')\Bigr\vert
\leqslant C_{\phi} \norme{x-x'}_\infty^{\alpha - \lfloor \alpha \rfloor}
\end{align}
for any $\beta \in \N^d$ such that $\beta = \lfloor \alpha\rfloor$ where $(\partial^\beta \cdot{} /\partial x^\beta)$ stands for the partial derivative of order $\beta$ and $\norme{\cdot{}}_{\infty}$ for the supremum norm on $[0,1]^d$.

\renewcommand{\theenumi}{\textbf{($\mathcal A$\arabic{enumi})}}
\renewcommand{\labelenumi}{\theenumi}
\begin{enumerate}
\setcounter{enumi}{3}
\item\label{hyp:lipschitz} The density $f_X $ of $X$ belongs to $ \mathcal C^{\alpha}([0,1]^d)$ for some $\alpha >0$.
\item\label{hyp:lipschitz_g} The regression function $m$ belongs to $\mathcal C^{\alpha}([0,1]^d)$.
\end{enumerate}
\renewcommand{\theenumi}{\alph{enumi}}
\renewcommand{\labelenumi}{\theenumi.}

Notice that, if the model $G$ belongs to ${\mathcal C}^{\alpha}$, then Assumption \ref{hyp:lipschitz_g} is satisfied.

\renewcommand{\theenumi}{\textbf{($\mathcal A$\arabic{enumi})}}
\renewcommand{\labelenumi}{\theenumi}
\begin{enumerate}
\setcounter{enumi}{5}
\item\label{hyp:kernel} 
Let $k\colon  {[0,1]} \to \mathbb{R}$ be a univariate kernel such that $\norme{k}_{\infty}< \infty$ and $\int_0^1 k({v})d{v}=1$.
We assume that $k$ is of order $(\lfloor \alpha \rfloor +1)$ which means that $\int_0^1 {v}^\ell k({v})d{v}=0$ for any $\ell \in \N$ such that $0< \ell\leq \lfloor \alpha \rfloor$ and $\int_0^1 {v}^{\lfloor \alpha \rfloor +1} k({v})d{v}\neq 0$. 
  {Furthermore, we define the multivariate kernel $K\colon [0,1]^d\to \R$ as: $K({v})=\prod_{k=1}^d k({v}_k)$ for any ${v}=({v}_1,\dots, {v}_d)\in [0,1]^d$.} Finally, we define $K_h({v})=K({v}/h)/h^d=\prod_{k=1}^d k({v_k}/h)/h^d$ for any ${v}=({v}_1,\dots, {v}_d)\in [0,1]^d$.
  \end{enumerate}
\renewcommand{\theenumi}{\alph{enumi}}
\renewcommand{\labelenumi}{\theenumi.}

Observe that $\int_{[0,1]^d}{v}^\beta K({v})d{v}=0$ for any $\beta \in \N^d$ such that $0< |\beta|\leq \lfloor \alpha \rfloor$ and $\int_{[0,1]^d}{v}^\beta K({v})d{v}\neq 0$ for some $\beta$ such that $|\beta|=\lfloor \alpha \rfloor +1$.

\renewcommand{\theenumi}{\textbf{($\mathcal A$\arabic{enumi})}}
\renewcommand{\labelenumi}{\theenumi}
\begin{enumerate}
\setcounter{enumi}{6}
\item\label{hyp:bandwidth} The sequence $(h_n)_{n\in \N}$ of bandwidths is positive and such that $h_n \to 0$ as $n\to\infty$. 
\end{enumerate}
\renewcommand{\theenumi}{\alph{enumi}}
\renewcommand{\labelenumi}{\theenumi.}

{The previous assumptions are pretty standard. See for instance \cite{tsybakov2009nonparametric}. More precisely, Assumptions ($\mathcal A_1$), ($\mathcal A_2$), ($\mathcal A_4$), ($\mathcal A_5$), ($\mathcal A_6$), and ($\mathcal A_7$) are common in the literature on nonparametric kernel regression estimation.  Assumption $\mathcal A_3$ naturally stems in our context where we aim to prove a central limit theorem for a quadratic functional of the regression function.}\\

In the rest of the paper, $C$ denotes a generic constant (deterministic and finite) which may vary from line to line.

\section{Estimation procedure and preliminary results} \label{sec:estimation_procedure}

In this section, we propose an estimator of $T$ based on two main ingredients: (a) estimation based on the efficient influence function of $T$ (see, e.g., \cite{bickel1993efficient,doksum1995nonparametric,van2000asymptotic}) and (b) mirror-type kernel estimators (see, e.g., \cite{bertin2020adaptive,pujol2022nonparametric}).

Let us first exhibit the general form of the estimator of $T$ considered in this paper:
\begin{equation}\label{estgeneralb}
T_n=\frac{1}{n}\sum_{i=1}^n \left(2Y_i-m_n(X_i)\right)m_n(X_i).
\end{equation} 
Here, $m_n$ is an estimator of the regression function $m$ of $Y$ on $X$ (with properties to be discussed later), and $\left(X_i,Y_i\right)_{1 \leq i \leq n}$ are independent copies of $(X,Y)$.

Such a form can actually be justified and inferred from a heuristic based on efficient influence functions. Indeed, let $\mathcal{P}$ be the set of  probability distributions on $[0,1]^d \times \mathbb{R}$ satisfying Assumptions \ref{hyp:densityb}, \ref{hyp:bounds}, \ref{hyp:lipschitz} and \ref{hyp:lipschitz_g}. Then we know from \cite{doksum1995nonparametric} that
\begin{align*}
(x,y)\mapsto (2y-m(x))m(x)-\E[\E[Y\vert X]]
\end{align*}
is the efficient influence function of $T$ under the nonparametric model $\mathcal{P}$ (see \cite{klein2024note_efficiency} for detailed computations). Thus, if the probability distribution of $(X,Y)$ is in $\mathcal{P}$ and if $m$ is known, we can state, from \cite[Equation (25.22)]{van2000asymptotic}, that the estimator 
\begin{equation}\label{estoracleb}
T_n=\frac{1}{n}\sum_{i=1}^n \left(2Y_i-m(X_i)\right)m(X_i) 
\end{equation} 
is  asymptotically efficient with 
optimal variance equal to $\textup{Var}\left(\left(2Y-m(X)\right)m(X)\right)$.
Unfortunately, the estimator in \eqref{estoracleb} is only an oracle, as $m$ is unknown in our case and needs then to be estimated, but this explains the intuition behind our focus on \eqref{estgeneralb}. 
As for the choice of the estimator of $m$, since we assume that the domain of the inputs is compact, the crucial point is to handle possible boundary effects. To do so, in \cite{doksum1995nonparametric} the authors estimate a truncated version of $T$ defined as
$T^{\textup{trunc},\varepsilon}=\E[\E[Y|X]^2\ind_{X \in (\varepsilon,1-\varepsilon)^d}]$. Even if $T^{\textup{trunc},\varepsilon}\to T$ as $\varepsilon\to 0$ under mild assumptions, the practical tuning of the parameter $\varepsilon$ depends on the unknown function $G$ and its choice has a large impact, see Figure \ref{fig:DK} in the numerical experiment section.

Here, we therefore prefer the use of  plugging mirror-type kernel estimators for $m(\cdot)$ in \eqref{estoracleb}.
{Indeed, kernel-based estimators are known to induce a mass loss near the boundary. We call a mirror-type transformation any transformation that corrects this phenomenon by pushing back the mass in the interior of the domain.}
Mirror-type transformations have been used, e.g., in \cite{bertin2020adaptive,pujol2022nonparametric} in the framework of density estimation.
In the following, we adapt these works to propose two nonparametric estimators of $m$, denoted as $\widehat{m}$ and $\widetilde{m}$. We then prove that both estimators $\widehat{T}_n$ and $\tilde{T}_n$ satisfy a central limit theorem with optimal asymptotic variance $\textup{Var}\left(\left(2Y-m(X)\right)m(X)\right)$.

\subsection{A mirror-type kernel estimator for the regression function}\label{sec:estim_bertin}

To estimate the regression function $m$, we take inspiration from Naradaya-Watson estimator \cite{nadaraya1964estimating,watson1964smooth} and the kernel-based plug-in estimator
studied in \cite{gine2008simple}.
More precisely, we consider a leave-one-out kernel estimator together with a mirror-type transformation introduced in \cite{bertin2020adaptive} to bypass 
boundary issues.
More precisely, the regression function estimator
is defined as follows:
\begin{align}\label{eq:estim_regression_bertin}
\widehat m_{n,h_n,i}(X_i) = \frac{\sum_{j\neq i} Y_j K_{h_n}\circ A_{X_i}(X_j-X_i)}{\sum_{j\neq i}  K_{h_n}\circ A_{X_i}(X_j-X_i)} 
\ind_{\sum_{j\neq i}  K_{h_n}\circ A_{X_i}(X_j-X_i)\neq 0                 }
\end{align}
for a bandwidth $h_n>0$, a mirror-type transformation $A$, and a kernel $K_{h_n}$ satisfying Assumption \ref{hyp:kernel} presented in Section \ref{sec:nota}. Then, \eqref{estoracleb} is adapted as
\begin{align}\label{def:estim_bertin}
\widehat T_{n,h_n} = \frac 1n \sum_{i=1}^n (2Y_i-\widehat m_{n,h_n,i}(X_i))\widehat m_{n,h_n,i}(X_i)\, .
\end{align}
 As for the mirror-transformation,  for $x\in [0,1]^d$, one may consider for instance
\begin{align}\label{def:mirror}
A_x \colon \Bigl\{\begin{array}{ccc}
\R^d & \to &\R^d\\
{v=(v_1,\dots,v_d)} & \mapsto & (a_1(x_1){v_1}, \dots , a_d(x_d){v_d})
\end{array}
\end{align}
with $a_i(s)\defeq 1-2 \ind_{( \frac 12,1]}(s)\in \{-1,1\}$, see Figure \ref{fig:MB} for an illustration. {Observe that $\{A_x, x\in [0,1]^d\}$ is a finite subset of $GL_d(\R)$ with cardinality $\kappa = 2^d$ (where $GL_d(\R)$ is the general linear group on $\R$). Then we denote these $\kappa$ elements $A_1, \ldots, A_{\kappa}$.} Moreover, it satisfies   
\begin{enumerate}
\item[(i)] for any $\ell=1,\dots,\kappa$, $\abs{\det (A_{\ell})}=1$;
\item[(ii)] \textbf{Mirror property}: 
\begin{align}\label{hyp:mirror}
\forall x\in [0,1]^d,\;   x+A_x^{-1}([0,1/2]^d)\subset [0,1]^d.
\end{align} 
\end{enumerate}

The $\mathcal C^\alpha([0,1]^d)$ regularity assumption can be relaxed to piecewise regularity as soon as the number $k$ of pieces is finite. The cardinality of {$\{A_x, x\in [0,1]^d\}$} is then increased to $\kappa=(2k)^d$.

\begin{figure}[ht!]
\centering
\begin{tikzpicture}[scale=.45]\footnotesize
 \pgfmathsetmacro{\xone}{-.4}
 \pgfmathsetmacro{\xtwo}{ 16.4}
 \pgfmathsetmacro{\yone}{-.4}
 \pgfmathsetmacro{\ytwo}{16.4}

\begin{scope}<+->;
 \draw[black] (8,8) node[anchor=north east] {$0$};

 \draw[black] (15,8) node[anchor=south west] {$1$};

 \draw[black] (8,15) node[anchor=south west] {$1$};

\fill[red]  (8,15) circle(1.2mm) ;
\fill[red]  (15,8) circle(1.2mm) ;

\draw[forestgreen] (9,13) node[anchor=south west] {$s$};
\draw[forestgreen] (9,2) node[anchor=south west] {$A_x(s)$};

\fill[forestgreen]  (9,13) circle(1.2mm) ;
\fill[forestgreen]  (9,2) circle(1.2mm) ;

\draw[forestgreen] (12,9) node[anchor=south west] {$v$};
\draw[forestgreen] (12,7) node[anchor=north west] {$A_x(v)$};

\fill[forestgreen]  (12,9) circle(1.2mm) ;
\fill[forestgreen]  (12,7) circle(1.2mm) ;

\draw[blue] (14,14) node[anchor=north west] {$w$};
\draw[blue] (2,14) node[anchor=north west] {$A_y(w)$};

\fill[blue]  (14,14) circle(1.2mm) ;
\fill[blue]  (2,14) circle(1.2mm) ;

\draw[blue] (10,11) node[anchor=north east] {$z$};
\draw[blue] (6,11) node[anchor=north east] {$A_y(z)$};

\fill[blue]  (10,11) circle(1.2mm) ;
\fill[blue]  (6,11) circle(1.2mm) ;

  \draw[black,thick,->] (\xone, 8) -- (\xtwo, 8) node[right] {$v_1$};
  \draw[black,thick,->] (8, \yone) -- (8, \ytwo) node[above] {$v_2$};
\end{scope}

\begin{scope}[thick,red]
  \foreach \x in {15}
    \draw (16-\x, 16-\x) rectangle (15,15);
\end{scope}
\end{tikzpicture}
\caption{Mirror-type transformation {defined in \eqref{def:mirror}}  with $d=2$, for $x=(1/3,3/4)$, and for $y=(2/3,1/5)$.}\label{fig:MB}
\end{figure}

Let $\widehat f_{n,h_n,i}$ be the leave-one-out estimator of the density function $f_X$ based on the 
$(n-1)$-sample $(X_1,\cdots, X_{i-1},\cdots, X_{i+1},\cdots,X_n)$: 
\begin{align}\label{eq:estim_density_bertin}
\widehat f_{n,h_n,i}(x)= \frac{1}{n-1}\sum_{j\neq i}  K_{h_n}\circ A_x(X_j-x). 
\end{align}
The following lemmas establish classical controls on the bias and on the variance of $\widehat f_{n,h_n,i}$ for all $i\in \{1,\cdots,d\}$.  The proofs are postponed to Appendix \ref{app:proof_bertin}.

\begin{lem}\label{lem:norme_sup_bertin}
Under Assumptions \ref{hyp:domain}, \ref{hyp:lipschitz}, and \ref{hyp:kernel}, for all $i\in \{1,\cdots,d\}$,
\begin{align}\label{eq:norme_sup_bertin}
\norme{\mathbb{E}\left[\widehat f_{n,h_n,i}\right]-f_X}_{\infty}=\norme{\mathbb{E}\left[\widehat f_{n,h_n,1}\right]-f_X}_{\infty} = O\Bigl(h_n^{\alpha}\Bigr)\cdot 
\end{align}
\end{lem}

To control the variance, one needs to define a supplementary assumption.

\renewcommand{\theenumi}{\textbf{($\mathcal A$\arabic{enumi})}}
\renewcommand{\labelenumi}{\theenumi}
\begin{enumerate}
\setcounter{enumi}{7}
\item\label{hyp:VCclass} Assume that the family of functions
\begin{align}\label{def:VCclass}
\mathcal F=\{K_{h,x}\colon y \in [0,1]^d \mapsto K_h(y-x)\in \R,\; h>0, x\in [0,1]^d\}
\end{align}
is a uniformly bounded Vapnik-Chervonenkis (VC)-class of functions, i.e.\ 
there exists positive numbers $A$, $B$, and $v$ such that, for all $K\in 
\mathcal F$, $\norme{K}_{\infty}<\infty$ and the $\varepsilon$-covering number $\mathcal N(\mathcal F,L^2(Q),\varepsilon)$ of $\mathcal F$ for the $L^2(Q)$-distance
satisfies
\[
\mathcal N(\mathcal F,L^2(Q),\varepsilon)\leqslant \Bigl(\frac{AB}{\varepsilon}\Bigr)^v
\]
for every probability measure $Q$ on $[0,1]^d$ and for every $\varepsilon \in 
(0,B)$.
\end{enumerate}
\renewcommand{\theenumi}{\alph{enumi}}
\renewcommand{\labelenumi}{\theenumi.}

{Assumption $(\mathcal A_8$) is a common hypothesis, called the uniformly bounded VC
class assumption, that is required to prove convergence of kernel density estimator (see, e.g., \cite{gine2001consistency}). Such an assumption is satisfied for commonly used kernels as Gaussian, Epanechnikov, Uniform, as mentioned in \cite{kim2019uniform}.}


\begin{lem}\label{lem:norme_sup_bertin_var}
Let $\delta_n\in (0,1)$.
Under Assumptions { \ref{hyp:domain}, \ref{hyp:kernel}, and \ref{hyp:VCclass},} there exists a constant $C>0$ such that we have, with probability $1-\delta_n$, 
	\begin{align}\label{eq:norme_sup_bertinvar}
		\norme{\widehat f_{n,h_n,i} - \mathbb{E}\left[\widehat f_{n,h_n,i}\right]}_{\infty}^2 \leq C \ \frac{\log(\frac{1}{h_n})+\log(\frac {2}{\delta_n})}{nh_n^d}
			\end{align}
			for all $i\in \{1,\cdots,d\}$.
\end{lem}

\begin{rem}
Choosing $h_n=n^{-1/(2\alpha+d)}$ and $\delta_n=1/n$ in Lemmas \ref{lem:norme_sup_bertin} and \ref{lem:norme_sup_bertin_var}, we get the following bound with probability $1-1/n$: 
$$
	\norme{\widehat f_{n,h_n,i}-f_X}_{\infty} \leq C \ \sqrt{\log(n)}n^{-\alpha/(2\alpha+d)} $$
 which corresponds to the optimal nonparametric rate up to the log factor. \end{rem}

Now, we can control the mean integrated squared error of $\widehat f_{n,h_n,i}$ together with its infimum. 

\begin{cor}\label{cor:cond_bertin}
Under Assumptions \ref{hyp:domain}, \ref{hyp:densityb}, \ref{hyp:lipschitz},  \ref{hyp:kernel}, and \ref{hyp:VCclass},
one has, for all $i\in \{1,\cdots,d\}$,
\begin{align}\label{eq:mise_f_bertin}
\E\bigl[\int_{[0,1]^d} (\widehat f_{n,h_n,i}(x)-f_X(x))^2 dx \bigr]= o(n^{-1/2})
\end{align}
and 
\begin{align}\label{eq:hat_f_bounded_bertin}
\frac{1}{\inf_{x \in [0,1]^d}\left|\widehat f_{n,h_n,i}(x)\right|}=O_{\P}(1)
\end{align}
	as soon as 
	$nh_n^{2d}\to \infty$ and $nh_n^{4\alpha}\to  0$ as $n\to \infty$. 
\end{cor}

Assuming $h_n=n^{-\gamma}$ with $\gamma>0$, the conditions $nh_n^{2d}{\to} \infty$ and $nh_n^{4\alpha}{\to} 0$ lead to $1/(4\alpha)<\gamma<1/(2d)$ and $\alpha > d/2$. The proof of Corollary \ref{cor:cond_bertin} is postponed to Appendix \ref{app:proof_bertin}.

\subsection{An alternative regression estimator}\label{sec:estim_pujol}

{This alternative estimator is based on the mirror-image kernel density estimator proposed in \cite{pujol2022nonparametric} (see also \cite{liu2012exponential}), introduced as a generalization to dimension $d\geq 2$ of the estimator introduced in \cite{hominal1979estimation,schuster1985incorporating} for the dimension $d=1$ and generalized in \cite{behnen1985rank} to the dimension $d=2$ (see also \cite{gijbels1990estimating}).}
More precisely, we consider the following transformations: for any $z\in [0,1]$, $m^{-1}(z)=-z, \quad m^{0}(z)=z, \quad \text{and} \quad m^{1}(z)=2-z$ and, for any $a\in \{-1,0,1\}^d$ and $x\in [0,1]^d$,  the $d$-dimensional vector $M^a(x)=(m^{a_1}(x_1),\cdots,m^{a_d}(x_d))$. Then, let
\begin{align}\label{eq:estim_density_pujol}
\widetilde{f}_{n,h_n,i}(x) = \ind_{[0,1]^d}(x) \frac{1}{(n-1)} \sum_{j \neq i} \sum_{a\in\{-1,0,1\}^d} \widetilde K_{h_n}\bigl(M^a(X_j)-x\bigr),
\end{align}
and
\begin{align}\label{eq:estim_num_pujol}
\widetilde{g}_{n,h_n,i}(x) = \ind_{[0,1]^d}(x) \frac{1}{(n-1)} \sum_{j\neq i} Y_j \sum_{a\in\{-1,0,1\}^d}  \widetilde K_{h_n} \bigl(M^a(X_j)-x\bigr).
\end{align}
Here, the bandwidth $h_n>0$ and the multivariate kernel $\widetilde K$ satisfy Assumption \ref{hyp:kernelbis} stated below.
Now we propose the following regression function estimator:
\begin{align}\label{eq:estim_regression_pujol}
\widetilde{m}_{n,h_n,i}(X_i) = \frac{\widetilde{g}_{n,h_n,i}(X_i)}{\widetilde{f}_{n,h_n,i}(X_i)}
\end{align}
if $\widetilde{f}_{n,h_n,i}(X_i)\neq 0$ and $0$ otherwise. The associated plug-in estimator then becomes:
\begin{equation}\label{defalter:estim_doksum}
	\widetilde{T}_{n,h_n} = \frac 1n \sum_{i=1}^n (2Y_i-\widetilde{m}_{n,h_n,i}(X_i))\widetilde{m}_{n,h_n,i}(X_i).
\end{equation}

The assumptions on $f_X$ and $g$ are strengthened as follows.

\renewcommand{\theenumi}{\textbf{{($\mathcal A'$}\arabic{enumi})}}
\renewcommand{\labelenumi}{\theenumi}
\begin{enumerate}
	\setcounter{enumi}{3}
	\item\label{hyp:lipschitzbis} The density function $f_X \in \mathcal C^{\alpha}([0,1]^d)$ for some $\alpha >0$. Its derivatives of order $\beta$, with $0<\beta \leq \lfloor \alpha \rfloor$, vanish near the boundary.
 	 \end{enumerate}
{Assumption ($\mathcal A'$\arabic{enumi}) appears to be strict, and is satisfied for example by uniform and beta distributions. Nevertheless, from a practical point of view, we have performed several simulations (not reported here) using Gaussian distributions that violate this assumption, but with estimation accuracy similar to the uniform setting.}\\
\renewcommand{\theenumi}{\alph{enumi}}
\renewcommand{\labelenumi}{\theenumi.}

{We now list the assumptions on $K$ below.}

\renewcommand{\theenumi}{\textbf{{($\mathcal A'$}\arabic{enumi})}}
\renewcommand{\labelenumi}{\theenumi}
\begin{enumerate}
	\setcounter{enumi}{5}
	\item\label{hyp:kernelbis} 
Let $\widetilde k\colon  {[-1,1]} \to \mathbb{R}$ be a univariate kernel such that $\norme{\widetilde k}_{\infty}< \infty$ and $\int_{-1}^1 \widetilde k({v})d{v}=1$.
We assume that $\widetilde k$ is of order $(\lfloor \alpha \rfloor +1)$ which means that $\int_{-1}^1 {v}^\ell \widetilde k({v})d{v}=0$ for any $\ell \in \N$ such that $0< \ell\leq \lfloor \alpha \rfloor$ and $\int_{-1}^1 {v}^{\lfloor \alpha \rfloor +1} \widetilde k({v})d{v}\neq 0$.
  Furthermore, we define the multivariate
  kernel $\widetilde K\colon [-1,1]^d \to \R$ as: $\widetilde K({v})=\prod_{k=1}^d \widetilde k({v}_k)$ for any ${v}=({v}_1,\ldots,{v}_d)\in [-1,1]^d$. Finally, we define $\widetilde K_h({v})=\widetilde K({v}/h)/h^d=\prod_{k=1}^d \widetilde k({v}_k/h)/h^d$ for any ${v}=({v}_1,\ldots,{v}_d)\in [-1,1]^d$.

	 \end{enumerate}
\renewcommand{\theenumi}{\alph{enumi}}
\renewcommand{\labelenumi}{\theenumi.}

Both numerator and denominator of the estimator defined in \eqref{eq:estim_regression_pujol} are the sum of $3^d$ terms; one corresponds to the
original data in the unit hypercube $[0, 1]^d$, and each of the remaining terms corresponds to reflected
data across one of the $0$-faces, $1$-faces, \ldots , $d-1$-faces of the unit hypercube as illustrated in Figure \ref{fig:MI} for $d=2$. {Note that Figure \ref{fig:MI} is inspired from illustrations in \cite{gijbels2015flexible,pujol2022nonparametric}.}

\medskip

Note that the function $\widetilde{f}_{n,h_n,i}$
is supported on $[0, 1]^d$ and $\int_{[0,1]^d}\widetilde{f}_{n,h_n,i}(x)dx=1$.  Moreover, it satisfies the two following lemmas. 

\begin{figure}
	\centering
		
		\begin{tikzpicture}[scale=.3]\footnotesize
			\pgfmathsetmacro{\xone}{-1}
			\pgfmathsetmacro{\xtwo}{ 24}
			\pgfmathsetmacro{\yone}{-1}
			\pgfmathsetmacro{\ytwo}{24}
			
			\draw[draw=red,very thick, fill=red!15] (8,8) rectangle (15,15);

			\draw[draw=orange!30,fill=orange!30] (7.5,8.5) rectangle (10.5,11.5);
			\draw[draw=orange!30,fill=orange!30] (7.5,4.5) rectangle (10.5,7.5);
			\draw[draw=orange!30,fill=orange!30] (7.5,18.5) rectangle (10.5,21.5);

			\draw[draw=orange!30,fill=orange!30] (5.5,4.5) rectangle (8.5,7.5);
			\draw[draw=orange!30,fill=orange!30] (5.5,8.5) rectangle (8.5,11.5);
			\draw[draw=orange!30,fill=orange!30] (5.5,18.5) rectangle (8.5,21.5);
			
			\draw[draw=orange!30,fill=orange!30] (19.5,4.5) rectangle (22.5,7.5);
			\draw[draw=orange!30,fill=orange!30] (19.5,8.5) rectangle (22.5,11.5);
			\draw[draw=orange!30,fill=orange!30] (19.5,18.5) rectangle (22.5,21.5);

			\draw[draw=orange!60,fill=orange!60] (7.5,8.5) rectangle (8.5,11.5);
			\draw[draw=orange!60,fill=orange!60] (7.5,4.5) rectangle (8.5,7.5);
			\draw[draw=orange!60,fill=orange!60] (7.5,18.5) rectangle (8.5,21.5);

			\draw[black] (8,8) node[anchor=north east] {$0$};
			\begin{scope}<+->;
				\draw[black] (8,8) node[anchor=north east] {$0$};
				\draw[black] (15,8) node[anchor=south west] {$1$};
				\draw[black] (22,8) node[anchor=south west] {$2$};
				\draw[black] (8,15) node[anchor=south west] {$1$};
				\draw[black] (8,22) node[anchor=south west] {$2$};
				\draw[black] (0,8) node[anchor=south west] {$-1$};
				\draw[black] (8,0) node[anchor=south west] {$-1$};

				
				\draw[forestgreen] (9,10) node[anchor=north west] {$z$};
				\fill[forestgreen]  (9,10) circle(1.2mm) ;
				\fill[orange]  (9,6) circle(1.2mm) ;
				\fill[orange]  (9,20) circle(1.2mm) ;
				\fill[orange]  (7,6) circle(1.2mm) ;
				\fill[orange]  (7,10) circle(1.2mm) ;
				\fill[orange]  (7,20) circle(1.2mm) ;
				\fill[orange]  (21,6) circle(1.2mm) ;
				\fill[orange]  (21,10) circle(1.2mm) ;
				\fill[orange]  (21,20) circle(1.2mm) ;

				\draw[black,thick,->] (\xone, 8) -- (\xtwo, 8); 
				\draw[black,thick,->] (8, \yone) -- (8, \ytwo); 
				\draw[black,thick] (0,0) -- (22, 0); 
				\draw[black,thick] (0,15) -- (22, 15);
				\draw[black,thick] (0,22) -- (22, 22);
				
				\draw[black,thick] (0,0) -- (0,22); 
				\draw[black,thick] (15,0) -- (15,22);
				\draw[black,thick] (22,0) -- (22, 22);   
			\end{scope}
			
			\draw[draw=red,very thick] (8,8) rectangle (15,15);
			
		\end{tikzpicture}
		\begin{tikzpicture}[scale=.3]\footnotesize
			\pgfmathsetmacro{\xone}{-1}
			\pgfmathsetmacro{\xtwo}{ 24}
			\pgfmathsetmacro{\yone}{-1}
			\pgfmathsetmacro{\ytwo}{24}
			
			\draw[draw=red,very thick, fill=red!15] (8,8) rectangle (15,15);

			\draw[draw=forestgreen!80,fill=forestgreen!80] (8,8.5) rectangle (8.5,11.5);
			\draw[draw=forestgreen!30,fill=forestgreen!30] (8.5,8.5) rectangle (10.5,11.5);
			
			\draw[black] (8,8) node[anchor=north east] {$0$};
			\begin{scope}<+->;
				\draw[black] (8,8) node[anchor=north east] {$0$};
				\draw[black] (15,8) node[anchor=south west] {$1$};
				\draw[black] (22,8) node[anchor=south west] {$2$};
				\draw[black] (8,15) node[anchor=south west] {$1$};
				\draw[black] (8,22) node[anchor=south west] {$2$};
				\draw[black] (0,8) node[anchor=south west] {$-1$};
				\draw[black] (8,0) node[anchor=south west] {$-1$};

				
				\draw[forestgreen] (9,10) node[anchor=north west] {$z$};
				\fill[forestgreen]  (9,10) circle(1.2mm) ;

				\draw[black,thick,->] (\xone, 8) -- (\xtwo, 8); 
				\draw[black,thick,->] (8, \yone) -- (8, \ytwo); 
				\draw[black,thick] (0,0) -- (22, 0); 
				\draw[black,thick] (0,15) -- (22, 15);
				\draw[black,thick] (0,22) -- (22, 22);
				
				\draw[black,thick] (0,0) -- (0,22); 
				\draw[black,thick] (15,0) -- (15,22);
				\draw[black,thick] (22,0) -- (22, 22);   
			\end{scope}
			
			\draw[draw=red,very thick] (8,8) rectangle (15,15);
			
		\end{tikzpicture}
	\caption{Mirror-image transformation with $d=2$ in red. In the left-hand side of the figure, a data-point in $[0,1]^d$ (in green) and its 8 mirror-images (in orange). A kernel is fitted over all points of the
	augmented dataset (orange areas). The darker orange regions correspond to the regions where several kernels overlap. In the right-hand side of the figure, the data-point (in green) and the restriction to $[0,1]^d$ of the augmented dataset (in green).}\label{fig:MI}
\end{figure}

\begin{lem}  
	\label{lem:norme_sup_pujol}
	Let $h_n \in (0,1/2)$.
Under Assumptions the \ref{hyp:domain}, \ref{hyp:lipschitzbis}, and \ref{hyp:kernelbis},
	\begin{align}\label{eq:norme_sup_pujol}
\forall i\in \{1,\cdots,d\},\quad		\norme{\mathbb{E}\left[\widetilde{f}_{n,h_n,i}\right]-f_X}_{\infty}=\norme{\mathbb{E}\left[\widetilde{f}_{n,h_n,1}\right]-f_X}_{\infty} = O(h_n^{\alpha})\cdot \end{align}
\end{lem}

	Lemma \ref{lem:norme_sup_pujol} is an extension 
	of 	\cite[Proposition 5.2]{pujol2022nonparametric} to every $\alpha >0$ in Assumption \ref{hyp:lipschitzbis} (itself an extension 
	to every dimension and every $\alpha \in (0, 2]$ of \cite[Lemma 3.1]{liu2012exponential}). 
		Its proof is postponed to Appendix \ref{app:proof_pujol}.

\begin{lem} \cite[Proposition 5.3]{pujol2022nonparametric} 
	\label{lem:norme_sup_pujol_var}
Let $h_n \in (0,1/2)$ and  $\delta_n\in (0,1)$.
	Under the Assumptions { \ref{hyp:domain}, \ref{hyp:kernelbis}, and \ref{hyp:VCclass} and with the kernel $\widetilde K$,} there exists a constant $C>0$ such that we have, with probability $1-\delta_n$, 
\begin{align}\label{eq:norme_sup_bertinpuj}
\forall i\in \{1,\cdots,d\},\quad	\norme{\widetilde{f}_{n,h_n,i}-\mathbb{E}\left[\widetilde{f}_{n,h_n,i}\right]}_{\infty}^2 \leq C \ \frac{\left(\log(\frac{1}{h_n})\right)_+ +\log(\frac{2}{\delta_n})}{nh_n^d}.
\end{align}
\end{lem}

Once more, choosing $h_n=n^{-1/(2\alpha+d)}$ and $\delta_n=1/n$ in Lemmas \ref{lem:norme_sup_pujol} and \ref{lem:norme_sup_pujol_var}, we recover the optimal nonparametric rate up to the log factor. 
    As before, we can deduce from Lemmas \ref{lem:norme_sup_pujol} and \ref{lem:norme_sup_pujol_var} the control of the mean integrated squared error of $\widetilde{f}_{n,h_n,i}$ and its infimum.

\begin{cor}\label{lem:mise_pujol}
Assume that Assumptions \ref{hyp:domain}, \ref{hyp:densityb}, \ref{hyp:lipschitzbis}, and \ref{hyp:kernelbis} are satisfied, as well as Assumption \ref{hyp:VCclass} with kernel $\widetilde K$. Then one has, for all $i\in \{1,\cdots,d\}$,
\begin{align}\label{eq:mise_f_pujol}
	\E\Bigl[\int_{[0,1]d} (\widetilde{f}_{n,h_n,i}(x)-f_X(x))^2 dx \Bigr]= o(n^{-1/2}) \quad \textrm{and} \quad \frac{1}{\inf_{x \in [0,1]^d}\left|\widetilde{f}_{n,h_n,i}(x)\right|}=O_{\P}(1)
\end{align}
$nh_n^{2d}\to \infty$ and $nh_n^{4\alpha}\to 0$ as $n\to \infty$. 
\end{cor}

	Assuming $h_n=n^{-\gamma}$ with $\gamma>0$, the conditions $nh_n^{2d}{\to} \infty$ and $nh_n^{4\alpha}{\to} 0$ lead once more to $1/(4\alpha)<\gamma<1/(2d)$ and $\alpha > d/2$.

\medskip

	The proof of Corollary \ref{lem:mise_pujol} comes from Lemmas \ref{lem:norme_sup_pujol} and \ref{lem:norme_sup_pujol_var} and is similar to the proof of Corollary \ref{cor:cond_bertin} from  Lemmas \ref{lem:norme_sup_bertin} and \ref{lem:norme_sup_bertin_var}. It is thus skipped.

\section{Central limit theorems}\label{sec:main_compact}

Here, we prove a central limit theorem for both estimators defined in \eqref{def:estim_bertin} and \eqref{defalter:estim_doksum}.

\begin{thm}[Central limit theorem] \label{th:mainresult_compact}
(i) Under Assumptions \ref{hyp:domain} to \ref{hyp:VCclass}, one has
	\begin{align}\label{eq:mainresult_gauss_compact}
	\sqrt{n}\bigl (\widehat T_{n,h_n}-\E[\E[Y|X]^2]\bigr ) \xrightarrow[n \rightarrow  \infty]{\mathcal L} \mathcal{N}\bigl(0,\Var((2Y-m(X))m(X))\bigr)
	\end{align}
as soon as $\alpha > d/2$ and $h_n=n^{-\gamma}$ with $1/(4\alpha)<\gamma<1/(2d)$. 

(ii) Replacing Assumptions \ref{hyp:lipschitz} and \ref{hyp:kernel} respectively by Assumptions \ref{hyp:lipschitzbis} and \ref{hyp:kernelbis} {in (i)}, the same conclusion holds true for $\widetilde T_{n,h_n}$. \end{thm}

The proof of Theorem \ref{th:mainresult_compact} is postponed to Appendix \ref{app:mainresult_compact} while the proofs of the following results are postponed to Appendix \ref{app:maincor_compact}. 

\begin{prop}[Asymptotic efficiency for $\widehat T_{n,h_n}$ and $\widetilde T_{n,h_n}$]\label{prop:asymp_eff}
Under the assumptions of Theorems \ref{th:mainresult_compact}, {$\widehat{T}_{n,h_n}$} and $\widetilde{T}_{n,h_n}$ are asymptotically efficient to estimate $\E[\E[Y\vert X]^2]$ from an i.i.d.\ sample $(X_i,Y_i)_{i=1,\cdots,n}$ of the pair $(X,Y)$.
\end{prop}

Using the delta method \cite[Theorem 3.1]{van2000asymptotic}, we are now able to get the asymptotic behavior of the estimators of $S^X$ derived respectively from \eqref{def:estim_bertin} and \eqref{defalter:estim_doksum}. Let
\[
\widehat S_{n,h_n}\defeq \frac{\widehat T_{n,h_n}-\Bigl (\frac 1n \sum_{j=1}^n Y_j\Bigr )^2}{\frac 1n \sum_{j=1}^n Y_j^2-\Bigl (\frac 1n \sum_{j=1}^n Y_j\Bigr )^2} \quad \text{and} \quad \widetilde S_{n,h_n}\defeq \frac{\widetilde T_{n,h_n}-\Bigl (\frac 1n \sum_{j=1}^n Y_j\Bigr )^2}{\frac 1n \sum_{j=1}^n Y_j^2-\Bigl (\frac 1n \sum_{j=1}^n Y_j\Bigr )^2}.
\]
\begin{cor}[Central limit theorem and asymptotic efficiency for $\widehat S_{n,h_n}$ and $\widetilde S_{n,h_n}$] \label{cor:TCL_global}

\noindent\\
(i) Under the assumptions of Theorem \ref{th:mainresult_compact} for $\widehat T_{n,h_n}$, one has 
\begin{align}\label{eq:TCL_global}
\sqrt{n}\left(\widehat S_{n,h_n}-S^X \right) \xrightarrow[n \rightarrow  \infty]{\mathcal L} \mathcal{N}(0,\sigma^2),
\end{align}
where the limiting variance  $\sigma^2$ has an explicit expression {given in the proof in \eqref{eq:sigma2}}. 

(ii) Under the assumptions of Theorem \ref{th:mainresult_compact} for $\widetilde T_{n,h_n}$, the same result holds for $\widetilde S_{n,h_n}$.

(iii) Moreover, $\widehat S_{n,h_n}$ and $\widetilde S_{n,h_n}$ are asymptotically efficient to estimate $S^X$ from an i.i.d.\ sample $(X_i,Y_i)_{i=1,\cdots,n}$ of the pair $(X,Y)$.
\end{cor}

Using once more the delta method, we deduce the asymptotic behavior of the vector of the $p$ first-order Sobol' indices. Let us denote by $S^{i}$ the first-order Sobol' index associated to $X=V_i$ and its estimators $\widehat S^i$ and $\widetilde S^i$ given by:
\begin{align*}
\widehat S^i_{n,h_n} \defeq \frac{\widehat T_{n,h_n}-\Bigl (\frac 1n \sum_{j=1}^n Y_j\Bigr )^2}{\frac 1n \sum_{j=1}^n Y_j^2-\Bigl (\frac 1n \sum_{j=1}^n Y_j\Bigr )^2} \, , \; \widetilde S^i_{n,h_n} \defeq \frac{\widetilde{T}_{n,h_n}-\Bigl (\frac 1n \sum_{j=1}^n Y_j\Bigr )^2}{\frac 1n \sum_{j=1}^n Y_j^2-\Bigl (\frac 1n \sum_{j=1}^n Y_j\Bigr )^2} \, \cdot
\end{align*} 
{Under the assumptions of Theorem \ref{th:mainresult_compact}, one may prove that $(\widehat S^1_{n,h_n},\dots, \widehat S^p_{n,h_n})$ and \\
 $(\widetilde S^1_{n,h_n},\dots, \widetilde S^p_{n,h_n})$ satisfy a central limit theorem with limiting variances that can be computed explicitly performing classical computations in statistics. Being useless in the rest of the paper since we do not perform statistical testing, such explicit expressions are not provided here.  
Moreover, $(\widehat S^1_{n,h_n},\dots, \widehat S^p_{n,h_n})$ and $(\widetilde S^1_{n,h_n},\dots, \widetilde S^p_{n,h_n})$ are asymptotically efficient to estimate $(S^1,\dots,S^p)$ from an i.i.d.\ sample $(X_i,Y_i)_{i=1,\cdots,n}$ of the pair $(X,Y)$.}

\medskip

{The proof follows from \cite[Theorem 25.50]{van2000asymptotic} and  does not rely on the fact that we are dealing with first-order Sobol' indices. Hence, if $u_1,\ldots,u_r$ are distinct subsets of $\{1,\ldots,p\}$, we also have under the same assumptions that $(\widehat S^{u_1}_{n,h_n},\dots, \widehat S^{u_r}_{n,h_n})$ and $(\widetilde S^{u_1}_{n,h_n},\dots, \widetilde S^{u_r}_{n,h_n})$ satisfy a central limit theorem with limiting variances that can be computed explicitly. }

 {
\begin{rem}
Assumptions \ref{hyp:lipschitz}, \ref{hyp:lipschitzbis} and \ref{hyp:lipschitz_g} are isotropic regularity assumptions. The results of this paper could be easily extended to anisotropic classes of regularity at the cost of more complicated notation. Then the bandwidth $h_n$ in Assumptions \ref{hyp:kernel}, \ref{hyp:kernelbis}, \ref{hyp:bandwidth} and \ref{hyp:VCclass} would be replaced by $d$ distinct bandwidths $h_{n,k}$, one for each dimension of $X$. In the numerical Section \ref{sec:num}, we therefore select one bandwidth per input dimension, using leave-one-out on the regression function.
\end{rem}}

\section{Comparison with existing methods}\label{sec:compar}

In this section, we compare several estimators of $\E[\E[Y|X]^2]$  proposed in the literature in terms of asymptotic variance in the central limit theorem. {All of them rely on the same given-data design of experiment as in the setting of our work, except the Pick-Freeze estimator. \\
Recall also that our estimator is asymptotically efficient and its limiting variance is given in Theorem \ref{th:mainresult_compact} by
\begin{align}\label{eq:limiting_variance_T}
\Var((2Y-m(X))m(X))
=\Var(m^2(X))+4\E[m^2(X)\sigma^2(X)].
\end{align}
}

\paragraph{A seminal kernel-based method}

As mentioned in Section \ref{sec:estimation_procedure}, the authors of \cite{doksum1995nonparametric} consider a truncated version of $T$ defined as $T^{\textup{trunc},\varepsilon}=\E[\E[Y|X]^2\ind_{X \in (\varepsilon,1-\varepsilon)^d}]$ to bypass boundary issues, with limiting variance
\begin{align*}
\Var~(\ind_{X \in (\varepsilon,1-\varepsilon)^d}m^2(X))
+4\E[\ind_{X \in (\varepsilon,1-\varepsilon)^d}m^2(X)\sigma^2(X)]
\end{align*}
As mentioned before, notice that the fine tuning of $\varepsilon$ is cumbersome. Moreover,  the estimated quantity is not exactly $T$ but only its pruned version. 

\paragraph{An alternative kernel-based method in dimension one}

An alternative procedure has been introduced in \cite[page 11]{da2008efficient}.
Anyway, note that the estimator $\widehat T_n$ defined in \cite[page 11]{da2008efficient} is not easily tractable in practice. More precisely, the initial $n$-sample is split into two samples of sizes $n_1=\lfloor n/\log n\rfloor$ and $n_2=n-n_1\approx n$. The first sample is used to estimate the joint density of $(X,Y)$ while the second one is used in an outer loop to estimate the integral term by Monte-Carlo. 
 {The sequence of estimators $(\widehat T_n)_{n \geq 1}$ is proved to be asymptotically efficient to estimate $\E[\E[Y |X]^2]$ and its limiting variance is naturally given by \eqref{eq:limiting_variance_T} as expected.} We refer the reader to \cite[Theorems 3.4 and 3.5]{da2008efficient} and \cite{GGKL20}.

\paragraph{Nearest neighbor-based method}

One may also compare our results to the estimation procedure proposed in \cite{devroye2018nearest}, based on nearest neighbors. Here again, the initial $2n$-sample is split into two samples of equal size $n$.  The first sample allows to estimate the regression function $m(\cdot{})=\E[Y|X=\cdot{}]$ using the first nearest neighbor of $x$ among the points of the first sample while the second sample is used as a plug-in estimator. They proved that their estimator $S_n$ is consistent for any dimension $d$ of $X$ and that $\sqrt n (S_n-\E[S_n])$ is asymptotically Gaussian. { Nevertheless, 
 the bias term is negligible only if $d\leqslant 3$ and in the setting where $f_X$ is Lipschitz continuous and bounded away from zero. Thus it may happen that the central limit theorem does not hold true for $\sqrt n (S_n-\E[\E[Y|X]^2])$ if $d\geqslant 4$.} 
The limiting variance obtained with this estimation procedure is
\begin{align*}
2\Var(m^2(X))+5\E[m^2(X)\sigma^2(X)]+2\E[\sigma^2(X)],
\end{align*}
the multiplicative factor 2 in front of $\Var(m^2(X))$ taking into account that we have considered two samples of size $n$. When $d=1$, one may also use the estimation procedure based on ranks introduced in \cite{Chatterjee2019} and studied in \cite{GGKL20} with asymptotic variance
\begin{align*}
\Var(m^2(X))+4\E[m^2(X)\sigma^2(X)]+\E[\sigma^2(X)].
\end{align*}
Note that this asymptotic variance is not the efficient one, but this methodology has been proven to perform numerically particularly well.

\paragraph{Pick-Freeze method}

The limiting variance involved in the central limit theorem 
of the Pick-Freeze estimation, based on the particular Pick-Freeze design, is given by
\begin{align*}
\Var(m^2(X))+2\E[m^2(X)\sigma^2(X)]+\E[\sigma^2(X)].
\end{align*}

{Although the Pick-Freeze estimator and ours are both asymptotically  efficient, they rely on 
different types of design of experiments. 
So in some sense, the limiting variances are not comparable. Nevertheless, for a fixed sample size, one can compare numerically the numerical accuracy  of these estimation methods through the variability of the estimators (see the boxplots in Section \ref{sec:num}).}

\section{Numerical applications}\label{sec:num}

In this section, we illustrate the practical performances of one of our estimators on {two} analytical test cases {and one realistic flood case} coming from the sensitivity analysis literature. Recall that both estimators are based on a high-order kernel supported on $[0,1]^d$ for the first estimator and on $[-1,1]^d$ for the second one. From a theoretical perspective both of them lead to equivalent convergence results, but the first one actually suffers from strong numerical instabilities as discussed in the additional file that can be found on the following link   \url{https://hal.science/hal-04052837}. This explains why we focus on the second one in what follows. For all test cases:
\begin{itemize}
\item {We focus on the estimation of first-order and total Sobol' indices. Note that the estimation of the total Sobol' index associated to input $V_i$ relies on the estimation of the closed Sobol' index of order $p-1$ associated to $X=V_{\sim i}$ (see Equation \eqref{def:sobol_total} in Section \ref{sec:model}).}
\item We compute our second mirror-type estimator (\ref{defalter:estim_doksum}) with an Epanechnikov kernel of order 2 and 4 (see \cite{hansen2005exact} for a definition), with the kernel bandwidth being optimized via leave-one-out on the regression function ("Kernel 2" and "Kernel 4"). 
\item We also consider concurrent estimators, namely the nearest-neighbor estimator of \cite{devroye2018nearest} ("NN") and the asymptotically efficient version of Pick-Freeze estimator studied in \cite{janon2012asymptotic} ("PF1") for first-order indices and for total indices, and also the replicated version of Pick-Freeze estimator proposed in  \cite{tissot2015randomized} ("PF2"), the rank estimator of \cite{GGKL20} ("Rank"), and the lag estimator of \cite{klein2024efficiency} ("Lag") for first-order indices.
\item For all estimators, we generate a standard $n$-sample $(X_1,Y_1),\cdots,(X_n,Y_n)$ except for the Pick-Freeze method since it relies on a structured design of experiments.
\item Each experiment is repeated {100} times with a number of model evaluations fixed to $n=500$ (then to $n=1000$). The reference value {that approximates the true value of the index} is obtained from a Pick-Freeze estimation with very large sample size ($d\times 5.10^4$).
\end{itemize}

\begin{figure}[h]
\centering
\begin{subfigure}[b]{0.45\textwidth}
    \centering
    \includegraphics[width=\textwidth]{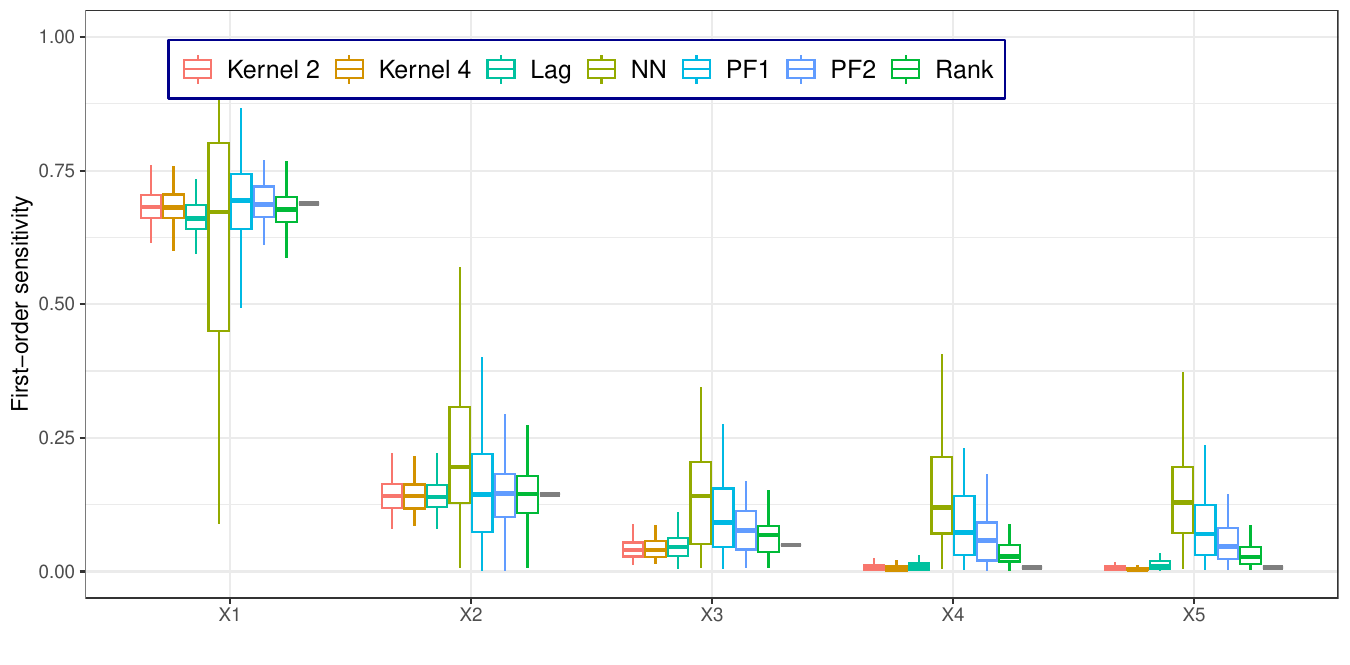}
    \caption{$n=500$}
\end{subfigure}
\hfill
\begin{subfigure}[b]{0.45\textwidth}
    \centering
    \includegraphics[width=\textwidth]{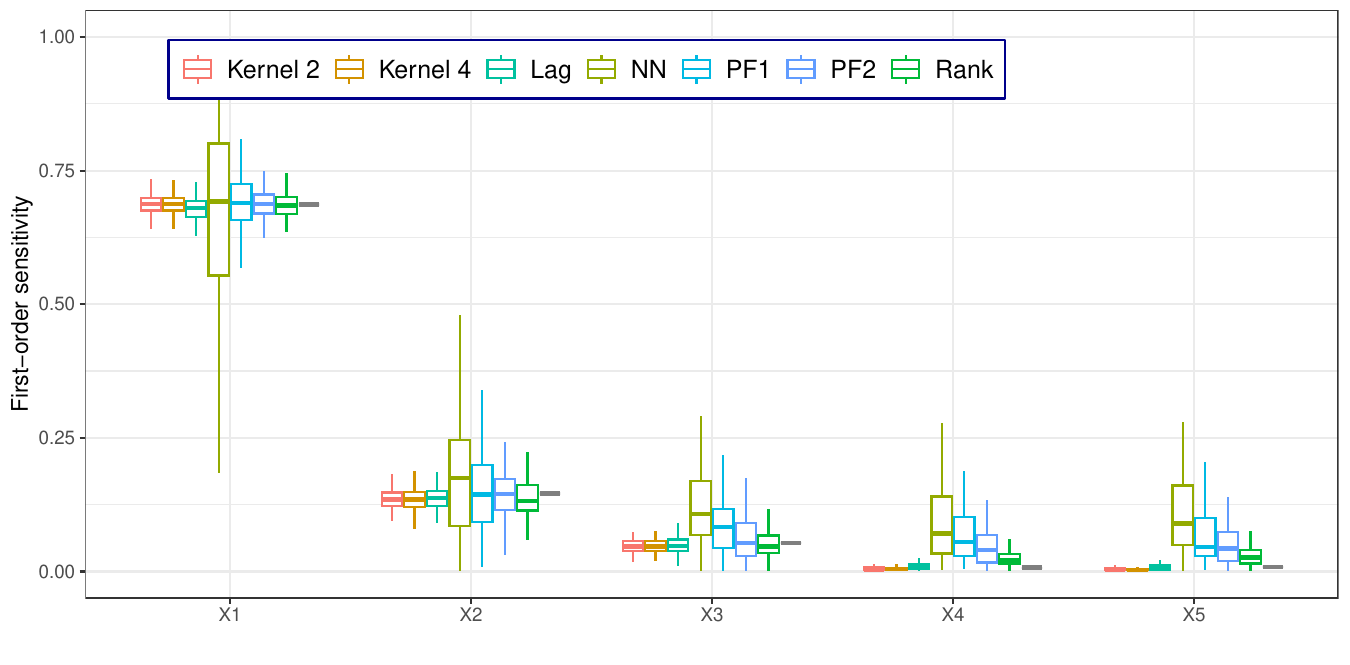}
    \caption{$n=1000$}
\end{subfigure}

\caption{Estimators for first-order indices of the Bratley function with $n=500$ (left) and $n=1000$ (right). The reference value is represented with a gray line.}
\label{fig:bratley1}
\end{figure}

\paragraph{The Bratley function}
First, we consider the Bratley function defined by:
\begin{equation}\label{eq:Bratley}
g_{\textrm{Bratley}}(V_1,\ldots,V_{p})=\sum_{i=1}^p (-1)^i \prod_{j=1}^{i} V_j,
\end{equation} 
with $V_i\sim \mathcal{U}([0,1])$ i.i.d. and $p=5$. The results of the $100$ experiments for all estimators of first-order indices are given as boxplots in Figure \ref{fig:bratley1}. We observe several trends.
\begin{itemize}
\item The nearest-neighbor estimator exhibits both large bias and large variance.
\item The Pick-Freeze estimators and the rank one perform well except when the sensitivity index is small.
\item Both the lag estimator and our mirror-type one have very small bias and variance.
\end{itemize}

\begin{figure}[h]
\centering
\begin{subfigure}[b]{0.45\textwidth}
    \centering
    \includegraphics[width=\textwidth]{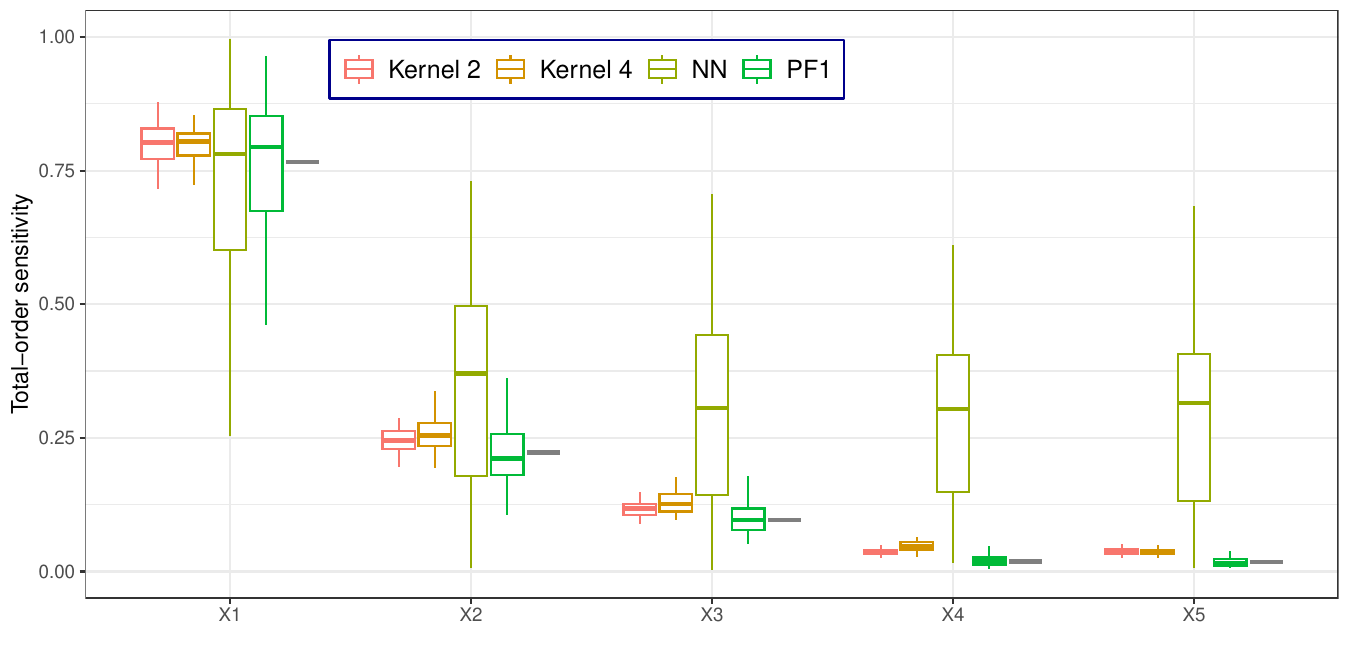}
    \caption{$n=500$}
\end{subfigure}
\hfill
\begin{subfigure}[b]{0.45\textwidth}
    \centering
    \includegraphics[width=\textwidth]{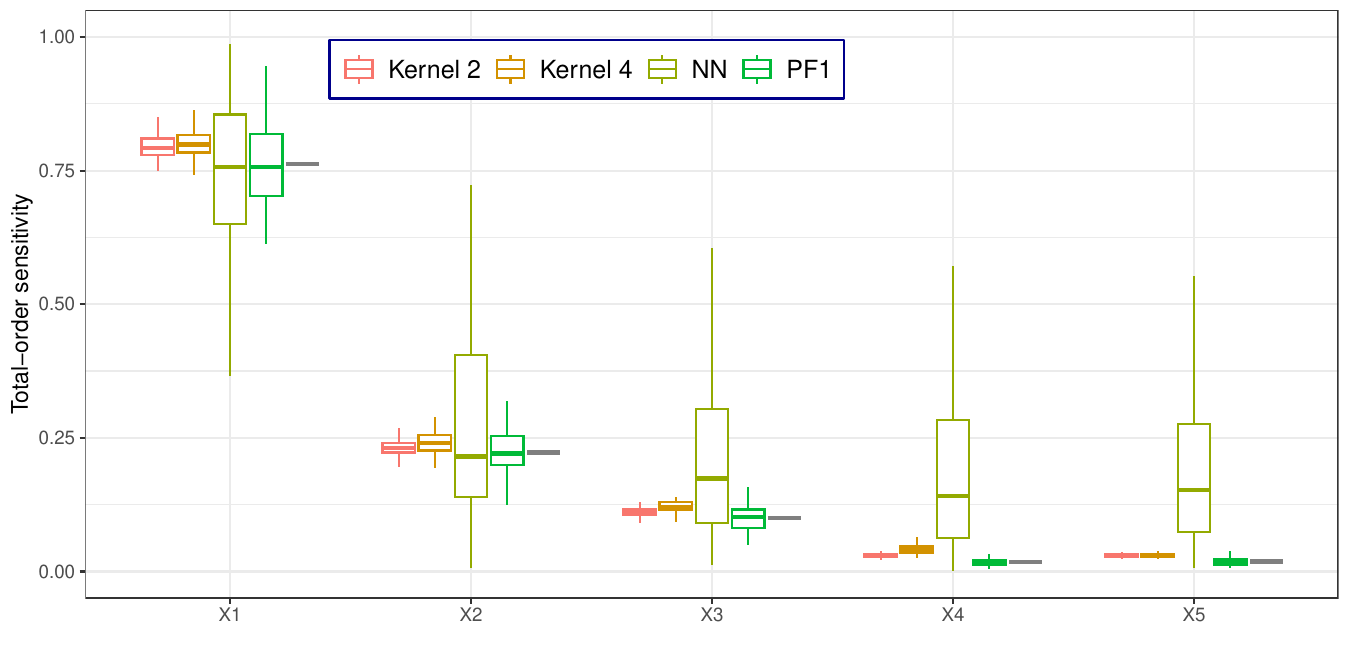}
    \caption{$n=1000$}
\end{subfigure}
\caption{Estimators for total indices of the Bratley function with $n=500$ (left) and $n=1000$ (right). The reference value of the index is represented with a gray line.}
\label{fig:bratleyT}
\end{figure}

For total indices in Figure \ref{fig:bratleyT}, recall that the only concurrent estimators are nearest-neighbor and Pick-Freeze estimators. Once again Pick-freeze estimators perform well, but the bias of the nearest-neighbor one is very large (here $d=p-1=4$), this bias increasing dramatically when the sensitivity index is small. Our mirror-type estimator still has small bias and very small variance for all input variables. Experiments with lower sample size are given in Appendix \ref{app:xps}.

\begin{figure}[h]
\centering
\begin{subfigure}[b]{0.45\textwidth}
    \centering
    \includegraphics[width=\textwidth]{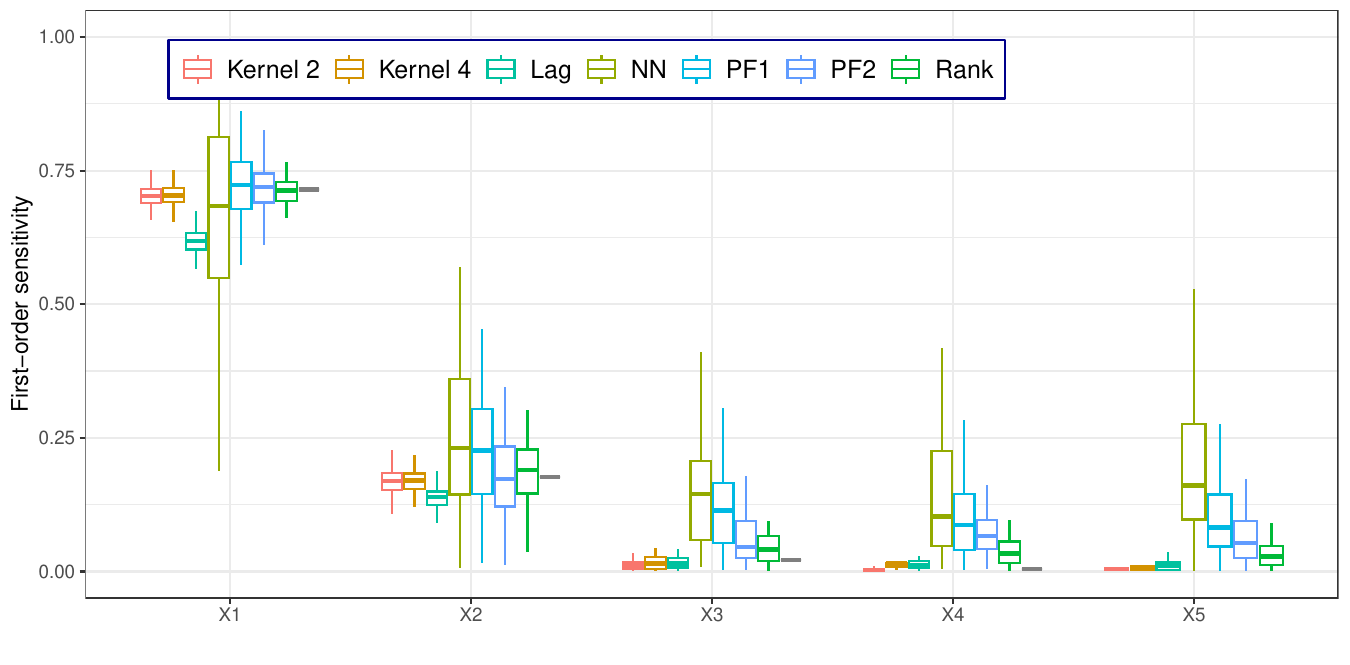}
    \caption{$n=500$}
\end{subfigure}
\hfill
\begin{subfigure}[b]{0.45\textwidth}
    \centering
    \includegraphics[width=\textwidth]{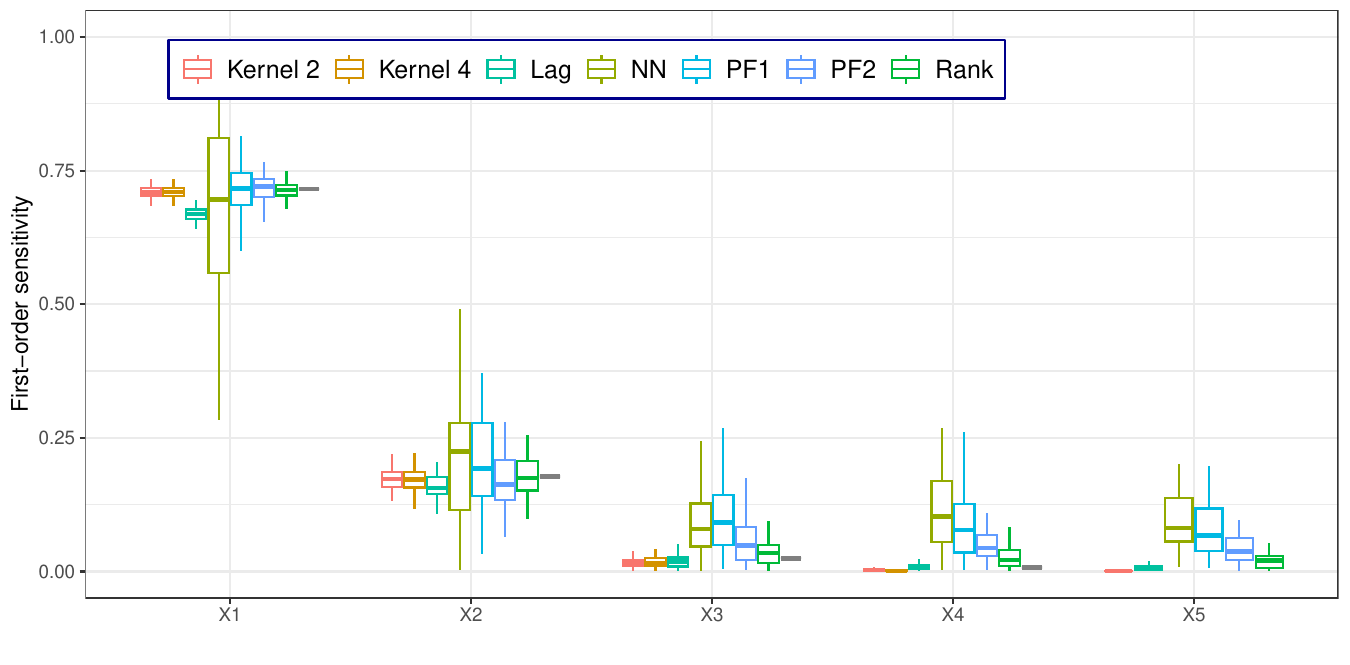}
    \caption{$n=1000$}
\end{subfigure}
\caption{Estimators for first-order indices of the g-Sobol function with $n=500$ (left) and $n=1000$ (right). The reference value of the index is represented with a gray line.}
\label{fig:gsobol1}
\end{figure}

\paragraph{The g-Sobol function}
Then, we investigate the g-Sobol function defined  by
\begin{equation}\label{eq:gsobol}
g_{\textrm{g-Sobol}}(V_1,\ldots,V_{p})= \prod_{i=1}^{p} \frac{\vert 4V_i-2\vert + a_i}{1+a_i},
\end{equation} 
with $V_i\sim \mathcal{U} ([0,1])$ i.i.d., $p=5$ and $a=(0,1,4.5,9,99)$. Notice that it is non-differentiable at any input value with a component equal to $0.5$, but the impact on our estimator performance is negligible, as can be seen in Figure \ref{fig:gsobol1} for first-order indices.
Except for the degraded performance of the lag estimator, the conclusions are the same as for the Bratley function, even for total indices displayed in Figure \ref{fig:gsobolT}. Experiments with lower sample size are given in Appendix \ref{app:xps}.

\begin{figure}[h]
\centering
\begin{subfigure}[b]{0.45\textwidth}
    \centering
    \includegraphics[width=\textwidth]{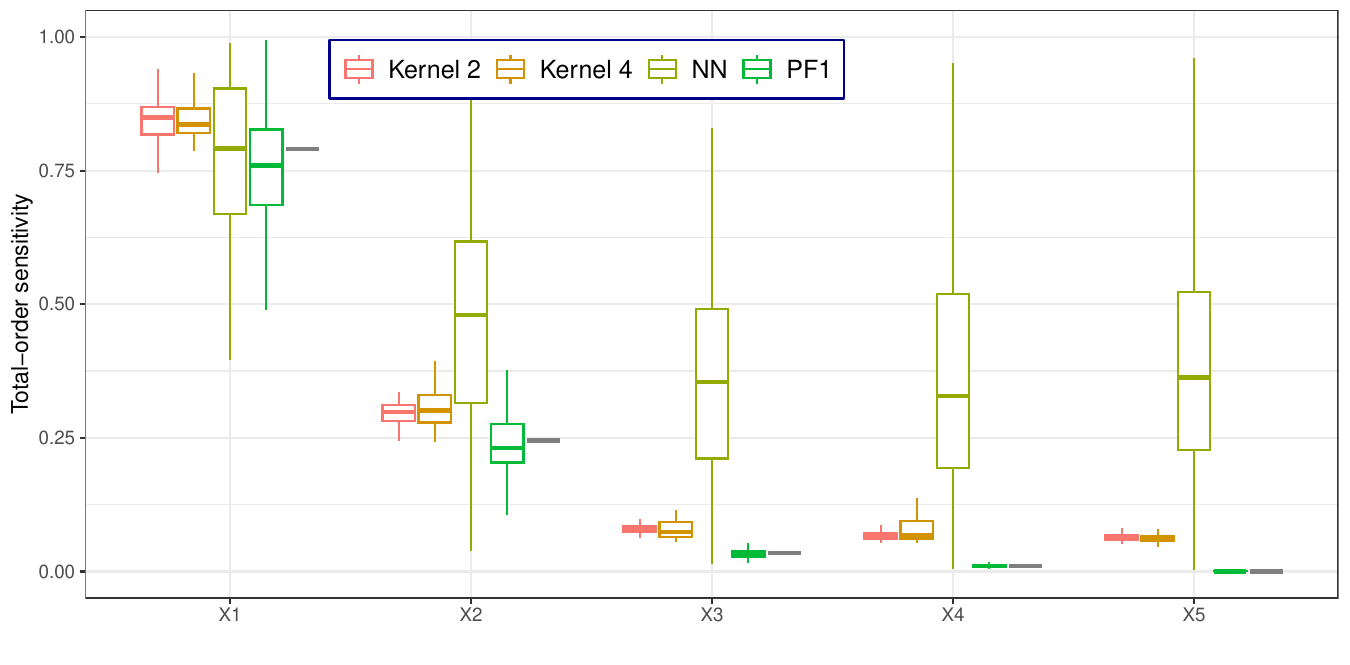}
    \caption{$n=500$}
\end{subfigure}
\hfill
\begin{subfigure}[b]{0.45\textwidth}
    \centering
    \includegraphics[width=\textwidth]{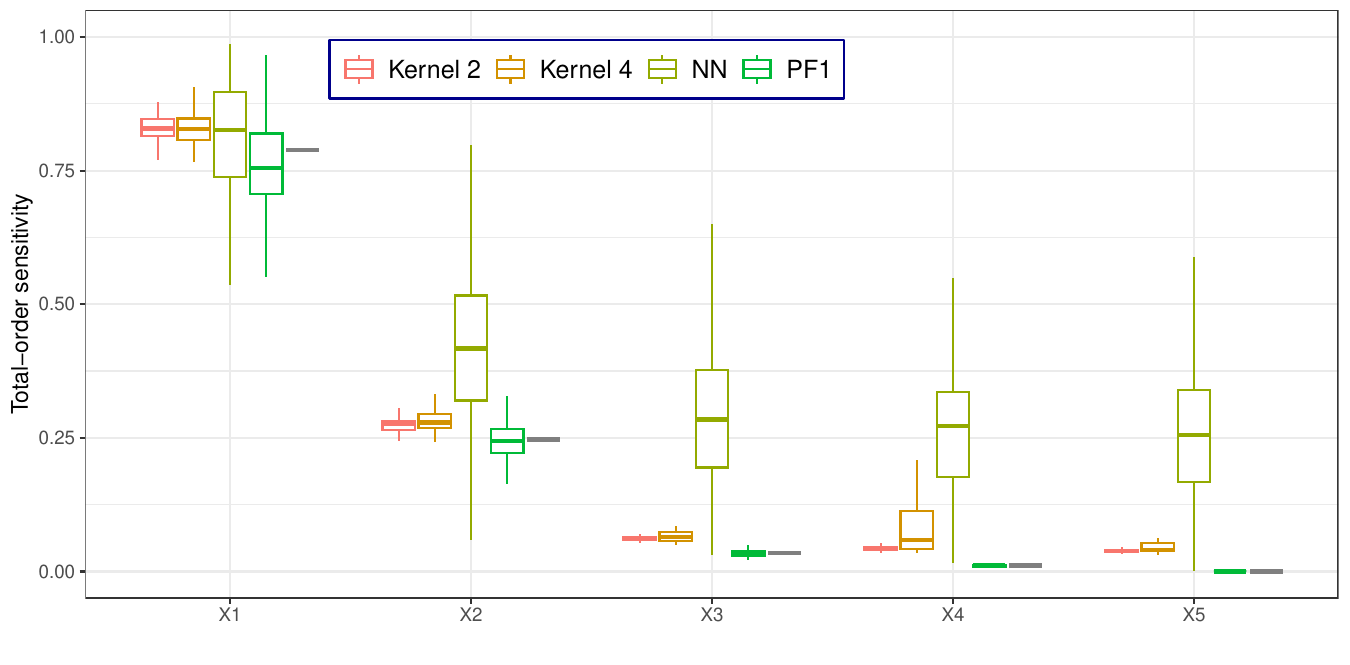}
    \caption{$n=1000$}
\end{subfigure}
\caption{Estimators for total indices of the g-Sobol function with $n=500$ (left) and $n=1000$ (right). The reference value of the index is represented with a gray line.}
\label{fig:gsobolT}
\end{figure}

{From Theorem \ref{th:mainresult_compact}, we know that it is necessary to choose a kernel of order $\lfloor \alpha \rfloor +1$ with $\alpha >d/2$ for the central limit theorem to hold. Here as $p=5$, the dimension $d$ of $X$ for total Sobol' indices is equal to $p-1=4$, thus the order of the kernel should be chosen at least equal to $3$. However we do not observe on Figure \ref{fig:gsobolT} better practical performance for a kernel of order 4 than for a one of order 2. Although the central limit theorem requires such a strong assumption about the order of the kernel, it is important to note that this assumption is not necessary for consistency results. Also we selected the bandwidth by cross-validation, that is with a data-driven approach, which naturally adapts to the order of the kernel. Currently, we have no clear idea of when it is beneficial to use a higher-order kernel in practice. It would be interesting to propose a strategy to choose adaptively and simultaneously the bandwidth and the order of the kernel as in \cite{bertin2020adaptive}. However this is out of the scope of the present work.}
\newpage

\paragraph{{A realistic flood model}}

{We now consider a realistic test case from the sensitivity analysis literature. The flood model used is a simplification of the 1D Saint-Venant hydrodynamic equations, assuming constant and uniform flows and very wide rectangular cross-sections. It consists of an equation involving the characteristics of the river section upstream of the industrial site:
\begin{equation*}
S = Z_v + H - H_d - C_b
\end{equation*}
where \( H \) is calculated as:
\begin{equation*}
H = \left(\frac{Q}{B K_s \sqrt{\frac{Z_m - Z_v}{L}}}\right)^{0.6}.
\end{equation*}
The model parameters are:
\begin{itemize}
    \item \( S \): overflow (in meters), model output;
    \item \( H \): maximum annual water level (in meters);
    \item \( Q \): maximum annual flow rate (in \( \mathrm{m}^3/\mathrm{s} \)), random input following a Gumbel max distribution \( Gu(1013, 558) \), truncated below at 500 and above at 3000;
    \item \( K_s \): strickler coefficient, random input following a normal distribution \(\mathcal N(30, 8) \), truncated below at 15;
    \item \( Z_v \): downstream riverbed elevation (in meters), random input following a triangular distribution \( T(49, 50, 51) \);
    \item \( Z_m \): upstream riverbed elevation (in meters), random input following a triangular distribution \( T(54, 55, 56) \);
    \item \( H_d \): dike height (in meters), random input following a uniform distribution \( \mathcal U(7, 9) \);
    \item \( C_b \): bank elevation (in meters), random input following a triangular distribution \( T(55, 55.5, 56) \);
    \item \( L \): length of the river section (in meters), random input following a triangular distribution \( T(4990, 5000, 5010) \);
    \item \( B \): river width (in meters), random input following a triangular distribution\\ \( T(295, 300, 305) \).
\end{itemize}
Among the model's input variables, \( H_d \) is considered random because it is studied as a design parameter. The other input variables are random due to their temporal and spatial variability, our lack of knowledge about their values, or imprecision in their estimation.} \\

\begin{figure}[h]
\centering
\begin{subfigure}[b]{0.45\textwidth}
    \centering
    \includegraphics[width=\textwidth]{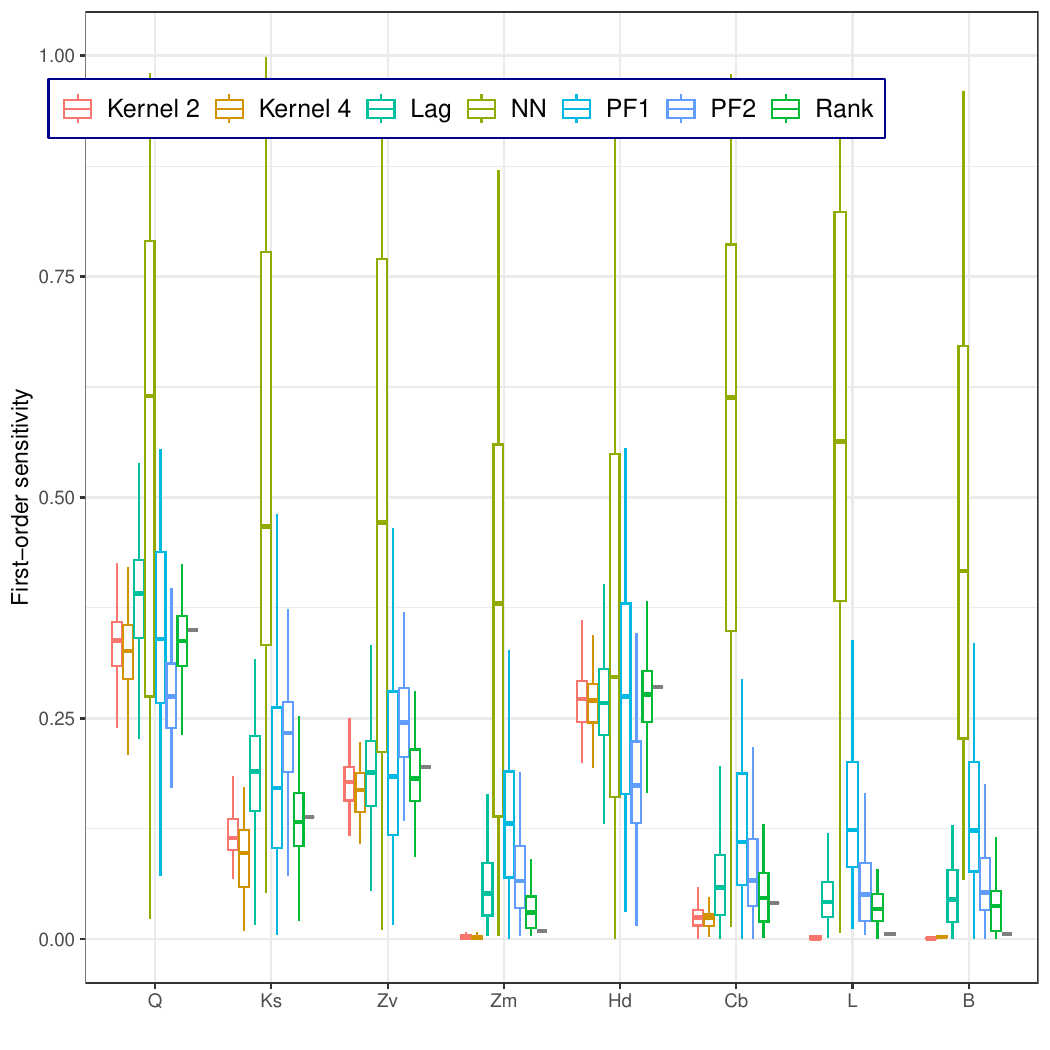}
    \caption{$n=500$}
\end{subfigure}
\hfill
\begin{subfigure}[b]{0.45\textwidth}
    \centering
    \includegraphics[width=\textwidth]{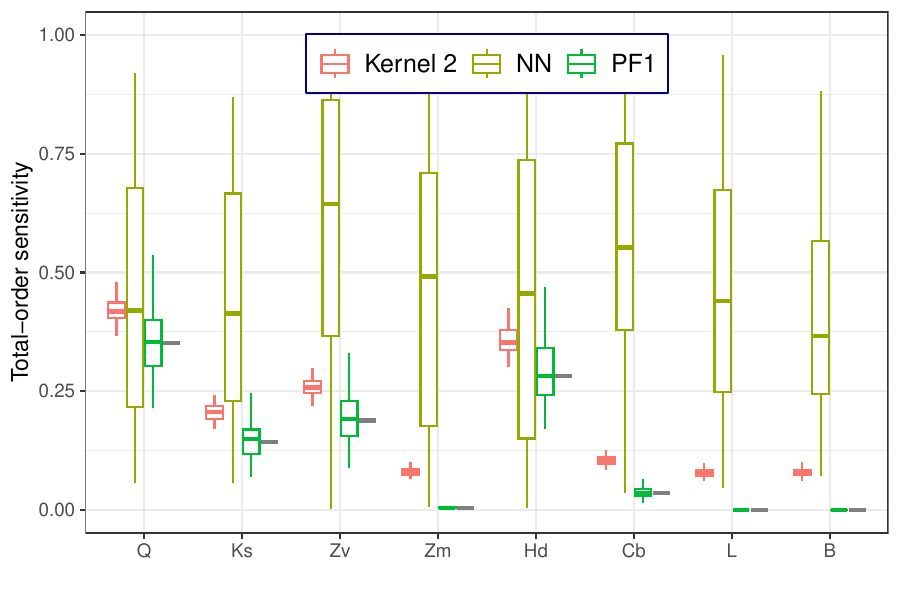}
    \caption{$n=500$}
\end{subfigure}
\caption{Estimators for first-order indices (left) and total indices (right) of the flood model with $n=500$. The reference value of the index is represented with a gray line.}
\label{fig:crue}
\end{figure}

{
We give in Figure \ref{fig:crue} the estimation of first-order and total Sobol' indices with a sample size $n=500$. Similarly to the previous analytical test cases, the nearest-neighbor estimator has very large bias and variance, and the Pick-Freeze one performs well. Concerning our kernel method, we present the results obtained with a kernel of order 2 or 4, and it clearly outperforms other estimators for first-order indices. For total ones, a residual bias persists, which should decrease by increasing the sample size. Nevertheless, the ranking of the inputs in terms of influence is preserved: the most important features are $Q$, $H_d$, and $Z_v$. These results are consistent with intuition, since it is expected that the overflow is mainly driven by the maximum annual flow rate $Q$ and the dike height $H_d$. As for the previous example, the results obtained for the estimation of total Sobol' indices with the kernel of order 4 are less favourable.  
}

\clearpage
\newpage
\mbox{~}

\appendix

\section{Proof of the results}\label{app:proof}

In the following, $h_n$ is simply denoted by $h$ and $C$ is a deterministic and finite constant, the value of which is allowed to change between occurrences. 
Recall also that $f_X$ is the density of $X$ with respect to the Lebesgue measure. The regression function is denoted by $m$: $m(x)=\E[Y\vert X=x]$ and we introduce the function $g$ defined by $g(x)=f_X(x)m(x)$. In addition, $\sigma^2(x)$ stands for $\Var(Y\vert X=x)$. 

\subsection{Proof of the preliminary results of Section \ref{sec:estim_bertin}}\label{app:proof_bertin}

Now, recall that $\widehat f_{n,h,i}$ defined in \eqref{eq:estim_density_bertin} and $\widehat m_{n,h,i}$ defined in \eqref{eq:estim_regression_bertin} are   respectively the leave-one-out estimator of $f_X$ based on $(X_1,\cdots, X_{i-1},\cdots, X_{i+1},\cdots,X_n)$ and the estimator of the regression function $m$ based on the 
$(n-1)$-input/output sample $((X_1,Y_1),\cdots,$ $(X_{i-1},Y_{i-1}),\cdots, (X_{i+1},Y_{i+1}),\cdots,(X_n,Y_n))$. Then, the leave-one-out estimator of $g$ is naturally given by $\widehat g_{n,h,i}=\widehat m_{n,h,i} \widehat f_{n,h,i}$.

\medskip

To lighten notation, we write $m_i$, $\sigma_i^2$, $g_i$, and $f_i$ for $m(X_i)$, $\sigma^2(X_i)$, $g(X_i)$, and $f_X(X_i)$ respectively. Additionally, 
$\widehat m_i$, $\widehat g_i$, and $\widehat f_i$ denote $\widehat m_{n,h,i}(X_i)$, $\widehat g_{n,h,i}(X_i)$, and $\widehat f_{n,h,i}(X_i)$ respectively. By abuse of notation, we denote $\widehat f_{n,h,i}(x)$ and $\widehat g_{n,h,i}(x)$ by $\widehat f_i(x)$ and $\widehat g_i(x)$.

\begin{proof}[Proof of 	Lemma \ref{lem:norme_sup_bertin}] The result can be deduced from the proof of the bound (8) in \cite[Proposition 1]{bertin2020adaptive} (see \cite[Section 7.2.]{bertin2020adaptive}) as the domain $[0,1]^d$ is compact. More precisely, 
	\begin{align*}
	 \E[\widehat f_{n,h,i}(x)]
		&= \int_{[0,1]^d} K_h\circ A_{x}(x'-x) f_X(x')dx'
=\int_{\mathcal D_x} K({v})   f_X(x+hA_x^{-1}({v}))  d{v}
	\end{align*}
where the last display is obtained after the variable change ${v}=A_{x}(x'-x)/h$, using the fact that $\abs{\det(A_x)}=1$ and where $\mathcal D_x{=A_x([0,1]^d-x)/h}$ is the integration domain after the variable change. Now observe that the mirror property in \eqref{hyp:mirror} ensures that, if $0<h \leqslant 1/2$ for all $i\in \{1,\cdots,d\}$, then $[0,1]^d  \subset {\mathcal D}_x$ for all  $x \in {[0,1]^d}$. Thus we have
	\begin{align}\label{eq:decomp1}
	\E[\widehat f_{n,h,i}(x)]-f_X(x)
		&= \int_{ [0,1]^d} K({v}) ( f_X(x+hA_x^{-1}({v}))-f_X(x))  d{v} 
\end{align}
since $\int_{[0,1]^d} K({v})d{v}=1$.  {If $\alpha \in (0,1)$, since $f_X\in \mathcal C^{\alpha}([0,1]^d)$, $\norme{K}_\infty<\infty$ and $\norme{A_x^{-1}({v})}_\infty\leqslant 1$, we conclude straightforwardly that
	\begin{align*}
	\abs{\E[\widehat f_{n,h,i}(x)]-f_X(x)}
		&\leqslant  C_{f_X} h^\alpha \int_{ [0,1]^d} \norme{A_x^{-1}({v})}_\infty^\alpha \abs{K({v})} d{v} \leqslant  C h^\alpha.
\end{align*} 
Now, for $\alpha\geqslant 1$,} observe that, still due to the mirror property in \eqref{hyp:mirror}, $0<h \leqslant 1/2$ for all $i\in \{1,\cdots,d\}$ ensures that $x+hA_x^{-1}({v}) \in [0,1]^d$ for all $x \in [0,1]^d$ and for all ${v}\in [0,1]^d$.
Thus the term in \eqref{eq:decomp1} can be handled with a Taylor expansion of $f_X$ (see, e.g., \cite[Theorem 5.4]{coleman2012calculus}). 
More precisely, by a multivariate Taylor expansion, since $f_X$ is $\mathcal C^{\alpha}([0,1]^d)$, we get: 
\begin{align*}
f_X(x+{v}h)-f_X(x)
&=\sum_{1\leqslant |\beta|< \lfloor \alpha \rfloor }  \frac{h^{|\beta|}}{\beta \, !} {v}^{\beta} \frac{\partial^{\beta}f_X}{\partial x^{\beta}}(x) + h^{\lfloor \alpha \rfloor} \sum_{|\beta|=\lfloor \alpha \rfloor}R_{\beta}(x+{v}h){v}^{\beta}
	\end{align*}
	with 
	\begin{align*}
	R_{\beta }(x+{v}h)=\frac {\lfloor \alpha \rfloor }{\beta !}\int _{0}^{1}(1-t)^{\lfloor \alpha \rfloor -1}\frac{\partial^{\beta}f_X}{\partial x^{\beta}}(x +t {v}h)\,dt.
	\end{align*}
Then,  {recalling that $(A_x^{-1}({v}))^\beta = (A_x^{-1}({v})_1)^{\beta_1}\dots (A_x^{-1}({v})_d)^{\beta_d}$}, we have
	\begin{align}
\E&[\widehat f_{n,h,i}(x)] -f_X(x)= \sum_{1\leqslant |\beta|< \lfloor \alpha \rfloor } \frac{h^{|\beta|}}{\beta \, !} \Bigl (\int_{[0,1]^d} (A_x^{-1}({v}))^\beta K({v}) d{v} \Bigr )  \frac{\partial^{\beta}f_X}{\partial x^{\beta}}(x) \notag\\
&\qquad\qquad\qquad\qquad\qquad  + h^{\lfloor \alpha \rfloor}\sum_{|\beta|=\lfloor \alpha \rfloor }    \int_{[0,1]^d} (A_x^{-1}({v}))^\beta K({v}) R_{\beta}(x+hA_x^{-1}({v})) d{v}\notag \\
	&= \lfloor \alpha \rfloor  h^{\lfloor \alpha \rfloor} \sum_{|\beta|=\lfloor \alpha \rfloor } \frac{1}{\beta !}   \int_{[0,1]^d} (A_x^{-1}({v}))^{\beta} K({v}) \int _{0}^{1}(1-t)^{\lfloor \alpha \rfloor-1}\frac{\partial^{\beta}f_X}{\partial x^{\beta}}(x +t hA_x^{-1}({v}))\,dt d{v}\label{eq:bias_kill} 
	\end{align}
	using the fact that $(A_x^{-1}({v}))^\beta$ is polynomial in $({v}_1,\dots,{v}_d)$ of degree $\beta$ and $K$ is of order $(\lfloor \alpha\rfloor +1)$. Now, using again that $K$ is of order $(\lfloor \alpha\rfloor +1)$, we get
	\begin{align*}
\E[\widehat f_{n,h,i}(x)] &-f_X(x)= \lfloor \alpha \rfloor  h^{\lfloor \alpha \rfloor} \sum_{|\beta|=\lfloor \alpha \rfloor } \frac{1}{\beta !} \int_{[0,1]^d} (A_x^{-1}({v}))^{\beta} K({v}) \\
	&\int _{0}^{1}(1-t)^{\lfloor \alpha \rfloor-1}\Bigl(\frac{\partial^{\beta}f_X}{\partial x^{\beta}}(x +t hA_x^{-1}({v})-\frac{\partial^{\beta}f_X}{\partial x^{\beta}}(x)\Bigr)\,dt d{v} .
	\end{align*}
Then,  using \eqref{eq:holder} since  $f_X\in\mathcal C^{\alpha}([0,1]^d)$, one gets 
	\begin{align*}
	\Bigl |\frac{\partial^\beta f_X}{\partial x^\beta}(x+thA_x^{-1}({v})) - \frac{\partial^\beta f_X}{\partial x^\beta}(x)\Bigr |
	\leqslant C_{f_X}( h t\|A_x^{-1}({v}) \|_{\infty})^{\alpha - \lfloor \alpha \rfloor}
	\end{align*}
	for all ${v} \in [0,1]^d$, $x\in [0,1]^d$, $t\in [0,1]$, $h\in (0,\infty)$, and $\beta \in \N^d$ such that $|\beta| = \lfloor \alpha\rfloor$. Then, 
	\begin{align*}
	 &|\E[\widehat f_{n,h,i}(x)] -f_X(x)| \\& \leqslant 
	C_{f_X} \lfloor \alpha \rfloor  h^{\alpha}\, \Bigl(\sum_{|\beta|=\lfloor \alpha \rfloor } \frac{1}{\beta !} \Bigr)   \Bigl(\int_{[0,1]^d} \|A_x^{-1}({v}) \|_{\infty}^\beta \|A_x^{-1}({v}) \|_{\infty}^{\alpha - \lfloor \alpha \rfloor} |K({v}) | d{v}\Bigr)\Bigl( \int _{0}^{1}(1-t)^{\lfloor \alpha \rfloor -1} t^{\alpha-\lfloor \alpha \rfloor } dt\Bigr)\notag\\
	&= { C  h^{\alpha}  \Bigl(\sum_{|\beta|=\lfloor \alpha \rfloor } \frac{1}{\beta !} \Bigr)  \Bigl(\int_{[0,1]^d}   |K({v}) | d{v}\Bigr) \Bigl( \int _{0}^{1}(1-t)^{\lfloor \alpha \rfloor -1} t^{\alpha-\lfloor \alpha \rfloor } dt\Bigr)
	\leqslant C h^{\alpha} }
	\end{align*}
	since   {$\|A_x^{-1}({v}) \|_{\infty}^\beta \leqslant 1$} and $\norme{K}_{\infty}<\infty$. 
\end{proof}

\begin{proof}[Proof of Lemma \ref{lem:norme_sup_bertin_var}]
For  {$u \subset \{1, \ldots , d \}$}, let $\sim u = \{1, \ldots , d \}\setminus u$.
For $x \in [0,1]^d$,
 \begin{align*}
&(n-1)\widehat f_{n,h,i}(x)=	\sum_{j\neq i}\prod_{\substack{1\leqslant k \leqslant d \text{~s.t.}\\ x_k \in [0,1/2]}} k_{h}
	(X_{k,j}-x_k)\prod_{\substack{1\leqslant k \leqslant d \text{~s.t.}\\ x_k \in (1/2,1]}} k_{h}(x_k-X_{k,j})\\\\
&	=  \sum_{j\neq i}\sum_{u \subseteq \{1, \ldots , d \}}
\prod_{k \in u} \ind_{[0,1/2]}(x_k)k_{h}(X_{k,j}-x_k) \prod_{k \in \sim u}\ind_{(1/2,1]}(x_k)k_{h}(x_k-X_{k,j})  \eqdef\sum_{j\neq i}K_{h,x}^*(X_j).
\end{align*} 
Arguing as in \cite[Section 5.1.2]{pujol2022nonparametric} and since the class of functions $\mathcal F$ defined by \eqref{def:VCclass} is a uniformly bounded VC-class of functions by \ref{hyp:VCclass}, we can say that the class $\mathcal F^*=\{K_{h,x}^*,\; h>0, x\in [0,1]^d\}$ is also a uniformly bounded VC-class of functions.
Thus it follows from a particular case of \cite[Corollary 13]{kim2019uniform} that the bound in \eqref{eq:norme_sup_bertinvar} holds with probability $1- \delta_n$.
\end{proof}

\begin{proof}[Proof of Corollary \ref{cor:cond_bertin}]
First, by Lemmas \ref{lem:norme_sup_bertin} and \ref{lem:norme_sup_bertin_var}, we get
\begin{align*}
\norme{\widehat f_{n,h,i}-f_X}_{\infty}^2&\leqslant 2\norme{\widehat f_{n,h,i} - \mathbb{E}\left[\widehat f_{n,h,i}\right]}_{\infty}^2 +2\norme{\mathbb{E}\left[\widehat f_{n,h,i}\right]-f_X}_{\infty}^2\\
&= O_{\P}\Bigl(h^{2\alpha} + \frac{\left(\log(\frac{1}{h_n})\right)_+ +\log(\frac{2}{\delta_n})}{nh_n^d} \Bigr)=o_{\P}(n^{-1/2})
\end{align*}
as soon as 
$nh^{2d}\to \infty$ and $nh^{4\alpha}\to 0$ by considering, e.g., $\delta_n=1/n$. 
Since the input space $[0,1]^d$ is compact by Assumption \ref{hyp:domain}, one concludes that \eqref{eq:mise_f_bertin} holds.

Moreover, for all $i\in \{1,\cdots,d\}$, as $\widehat f_{n,h,i}(x)=f_X(x)+\widehat f_{n,h,i}(x)-f_X(x)$, and from Assumption \ref{hyp:densityb}, we get with probability $1-1/\delta_n$ and for $n$ large enough:
\[
\inf_{x \in [0,1]^d}\left|\widehat f_{n,h,i}(x)\right| \geqslant \delta - C \ \Bigl(h^{2\alpha} + \frac{\left(\log(\frac{1}{h_n})\right)_+ +\log(\frac{2}{\delta_n})}{nh_n^d} \Bigr)
\]
	for some positive constant $C$. This last quantity is nonnegative as soon as $n$ is large enough when 
$nh^{d}\to \infty$ and $h^{\alpha}\to 0$ by considering, e.g., $\delta_n=1/n$. 
\end{proof}

We state Lemma \ref{lem:mise_g_bertin} below as a preliminary result for the proof of (i) of Theorem \ref{th:mainresult_compact}.

\begin{lem}\label{lem:mise_g_bertin}
Under Assumptions \ref{hyp:domain}, \ref{hyp:lipschitz_g},  \ref{hyp:kernel}, and \ref{hyp:VCclass},
\begin{align}\label{eq:mise_g_bertin}
\forall i\in\{1,\cdots,d\}, \quad \E\bigl[\int_{[0,1]^d} (\widehat g_{n,h,i}(x)-g(x))^2 dx \bigr]=o(n^{-1/2})
\end{align}
as soon as 
$nh^{2d}\to \infty$, and $nh^{4\alpha}\to 0$. 
\end{lem}

\begin{proof}[Proof of Lemma \ref{lem:mise_g_bertin}] Similarly as for density estimation, we first prove that $\norme{\mathbb{E}\left[\widehat g_{n,h,1}\right]-g}_{\infty} = O(h^{\alpha})$. Because $\E[\widehat g_{n,h,1}(x)] = \int_{[0,1]^d} K_h\circ A_{x}(x'-x) g(x')dx'$, we can follow the same lines as in the proof of \eqref{eq:norme_sup_bertin} since $g\in \mathcal C^{\alpha}([0,1]^d)$. Second, we prove that there exists some constant $C>0$ such that, with probability $1- \delta_n$,
	\begin{align}\label{eq:norme_sup_g_bertinvar}
	\norme{\widehat g_{n,h,1} - \mathbb{E}\left[\widehat g_{n,h,1}\right]}_{\infty} \leq C \ \Bigl(\frac{\log(\frac{1}{h})+\log(\frac{2}{\delta_n})}{nh^d}\Bigr).
	\end{align}
One has, for all $i\in \{1,\cdots,d\}$,	\begin{align*}
&(n-1)\widehat g_{n,h,i}(x)=	\sum_{j\neq i} Y_j\prod_{\substack{1\leqslant k \leqslant d\text{~s.t.}\\ x_k \in [0,1/2]}}k_{h}
	(X_{k,j}-x_k)\prod_{\substack{1\leqslant k \leqslant d \text{~s.t.}\\ x_k \in (1/2,1]}} k_{h}(x_k-X_{k,j})\\\\
&	=  \sum_{j\neq i}Y_j\sum_{u \subseteq \{1, \ldots , d \}}
\prod_{k \in u} \ind_{[0,1/2]}(x_k)k_{h}(X_{k,j}-x_k) \prod_{k \in \sim u}\ind_{(1/2,1]}(x_k)k_{h}(x_k-X_{k,j}) =O\Bigl( \sum_{j\neq i}K_{h,x}^*(X_j)\Bigr).
	\end{align*}
We conclude as in the proof of Lemma \ref{lem:norme_sup_bertin_var}. 
Finally, we deduce \eqref{eq:mise_g_bertin} as soon as 
$nh^{d}\to \infty$ and $h^{\alpha}\to 0$ by considering, e.g., $\delta_n=1/n$.
\end{proof}

\subsection{Proof of the preliminary results of Section \ref{sec:estim_pujol}}\label{app:proof_pujol}

Now, recall that $\widetilde f_{n,h,i}$ defined in \eqref{eq:estim_density_pujol} and $\widetilde g_{n,h,i}$ defined in \ref{eq:estim_num_pujol} are the leave-one-out estimators of $f_X$ and $g$ respectively based on the 
$(n-1)$-input/output sample $((X_1,Y_1),\cdots,$ $ (X_{i-1},Y_{i-1}),\cdots, (X_{i+1},Y_{i+1}),\cdots,(X_n,Y_n))$.
Then, the estimator of the regression function $m$ is naturally given by $\widetilde m_{n,h,i}=\widetilde g_{n,h,i}/\widetilde f_{n,h,i}$.

{ \begin{proof}[Proof of Lemma \ref{lem:norme_sup_pujol}]
It follows the same lines as in the proof of Lemma \ref{lem:norme_sup_bertin}. 
We consider the case $x \in [0,1/2]^d$. The remaining cases can be deduced by symmetry. Let ${\mathcal A} \defeq \{k \in \{1, \ldots , d\} \, : \, x_k \leq h\}$. For all $k \in \{1,\cdots,d\}$ and all $t_k \in [0,1]$, $k_{h}(M^{1}(t_k)-x_k)=\widetilde k_{h}(2-t_k-x_k)=0$ as $\widetilde k$ is supported on $[-1,1]$ and $h < 1/2$. In addition, for all $k \in \{1,\cdots,d\} \setminus \mathcal A$ and all $t_k \in [0,1]$,  $k_{h}(M^{-1}(t_k)-x_k)=\widetilde k_{h}(-t_k-x_k)=0$.

For any subset $\mathcal B \subset \mathcal A$, define $x_{\mathcal B}$ such that $x_{{\mathcal B},k}=x_k=M^0(x_k)$ if $k \notin \mathcal B$ and $x_{{\mathcal B},k}=-x_k=M^{-1}(x_k)$ if $k \in \mathcal B$. Then, the expected value of $\widetilde{f}_{n,h,i}(x)$ can be written as:
\begin{align*}
	\mathbb{E}[\widetilde{f}_{n,h,i}(x)]& = \sum_{\mathcal B \subset \mathcal A}\int_{[0,1]^d}\prod_{k \in \mathcal B}\widetilde k_{h}( -z_k-x_k) \prod_{k \notin \mathcal B}\widetilde k_{h}( z_k-x_k)f_X(z)dz \nonumber \\
& 	=  \sum_{\mathcal B \subset \mathcal A} \int_{{\mathcal X}_{\mathcal B}} \prod_{k =1}^{d} \widetilde k(u_k)f_X(z_{{\mathcal B}})du 
=   {\sum_{\mathcal B \subset \mathcal A} \int_{{\mathcal X}_{\mathcal B}} \widetilde K({v})f_X(z_{{\mathcal B}})d{v}}
	\end{align*}
with
\[
z_{{\mathcal B},k}=\begin{cases}
-x_k-{v}_kh, & \text{$k \in \mathcal B$}\\
x_k+{v}_kh, & \text{$k \notin \mathcal B$}
\end{cases} \text{~and~} {\mathcal X}_{{\mathcal B},k}=\begin{cases}
\{{v}_k \in [-1,1) \, : \, -x_k-{v}_kh \in [0,1]\}, & \text{$k \in \mathcal B$}\\
\{{v}_k \in [-1,1) \, : \, x_k+{v}_kh \in [0,1]\}, & \text{$k \notin \mathcal B$}
\end{cases}
\]
(recalling that the support of $\widetilde K$ is  $[-1,1]^d$) and ${\mathcal X}_{\mathcal B}=\prod_{k=1}^d {\mathcal X}_{{\mathcal B},k}$. Note that $({\mathcal X}_{\mathcal B})_{\mathcal B \subset \mathcal A}$ forms a partition of $[-1,1)^d$. Indeed, one has
\begin{align*}
	{\mathcal X}_{\mathcal B} = \prod_{k=1}^d {\mathcal X}_{{\mathcal B},k} = \prod_{k \in \mathcal B}[-1,-\frac{x_k}{h}) \, \prod_{k \in \mathcal A \setminus \mathcal B}[-\frac{x_k}{h},1)\, \prod_{k \in \{ 1,\cdots,d\}\setminus \mathcal A}[-1,1).
	\end{align*} 
	 {Let us now define $f_X^{\textup{MI}}$ on $[-1,2]^d$ such that, for all $y \in [-1,2]^d$, 
	\[
	f_X^{\textup{MI}}(y)=f_X(y^{\textup{MI}}) \quad \text{where} \quad y^{\textup{MI}}_k=\begin{cases}
	-y_k & \text{if $y_k\in [-1,0]$}\\
	y_k & \text{if $y_k\in [0,1]$}\\
	2-y_k & \text{if $y_k\in [1,2]$}\\
	\end{cases}.
	\]
Then, for all $x \in [0,1]^d$ and $a \in \{-1,0,1\}^d$, 
$f_X^{\textup{MI}}(M^a(x))=f_X(x)$ and similarly, for any $\mathcal B \subset \mathcal A$, $f_X^{\textup{MI}}(x_{\mathcal B})=f_X(x)$.
Now, since $f_X^{\textup{MI}}=f_X$ on $[0,1]^d$ and for any $w\in {\mathcal X}_{\mathcal B}$,  $f_X^{\textup{MI}}(z_{{\mathcal B}})=f_X(z_{{\mathcal B}})$, we have
\begin{align}\label{esp}
	\mathbb{E}[\widetilde{f}_{n,h,i}(x)]& = \sum_{\mathcal B \subset \mathcal A}\int_{{\mathcal X}_{\mathcal B}}
\widetilde K({v})f_X^{\textup{MI}}(z_{{\mathcal B}})d{v}
= \sum_{\mathcal B \subset \mathcal A}\int_{{\mathcal X}_{\mathcal B}'}
\widetilde K({w})f_X^{\textup{MI}}(z_{{\mathcal B}})d{w}
	\end{align}
using the variable change ${w}={v}_{\mathcal B}$, the symmetry of $\widetilde K$ and with  
\[
{\mathcal X}_{{\mathcal B},k}'=\begin{cases}
\{{w}_k \in [-1,1) \, : \, -x_k+{w}_kh \in [0,1]\}, & \text{$k \in \mathcal B$}\\
\{{w}_k \in [-1,1) \, : \, x_k+{w}_kh \in [0,1]\}, & \text{$k \notin \mathcal B$}
\end{cases}.
\]
 In addition, since $\int_{[-1,1]^d}\widetilde K({v})d{v}=1$, $({\mathcal X}_{\mathcal B}')_{\mathcal B \subset \mathcal A}$ forms also a partition of $[-1,1)^d$, and $f_X^{\textup{MI}}(x_{\mathcal B})=f_X(x)$, one has}
\begin{align}\label{truth}
	f_X(x)&= f_X(x) \int_{[-1,1]^d}  {\widetilde K({v})} d{v}  =
\sum_{\mathcal B \subset \mathcal A} \int_{{\mathcal X}_{\mathcal B}'}  {\widetilde K({v})} f_X(x)d{v} 
=\sum_{\mathcal B \subset \mathcal A} \int_{{\mathcal X}_{\mathcal B}'}  {\widetilde K({v})} f_X^{\textup{MI}}(x_{\mathcal B})d{v} \, .
	\end{align}
From \eqref{esp} and \eqref{truth}, we deduce
\[
\mathbb{E}[\widetilde{f}_{n,h,i}(x)]-f_X(x)=\sum_{\mathcal B \subset \mathcal A} \int_{{\mathcal X}_{\mathcal B}'}  {\widetilde K({v})} \Bigl(f_X^{\textup{MI}}(x_{\mathcal B}+vh)-f_X^{\textup{MI}}(x_{\mathcal B})\Bigr)d{w}. 
\]
By construction and since Assumption \ref{hyp:lipschitzbis} is satisfied (in particular derivatives up to order $\lfloor \alpha \rfloor$ vanish near the boundary), $f_X^{\textup{MI}}$ has derivatives up to order $\lfloor \alpha \rfloor$. 
Then, $f_X^{\textup{MI}}$ belongs to ${\mathcal C}^{\alpha}([-1,2]^d )$ with  {$C_{f_X^{\textup{MI}}}=3C_{f_X}$}. Thus, proceeding as in the proof of Lemma \ref{lem:norme_sup_bertin} we write:
\begin{align}\label{eq:taylor_fMI}
	f_X^{\textup{MI}}(x_{\mathcal B}+{v}h)-f_X^{\textup{MI}}(x_{\mathcal B})
	&=\sum_{1\leqslant |\beta|< \lfloor \alpha \rfloor }  \frac{h^{|\beta|}}{\beta \, !} {v}^{\beta} \frac{\partial^{\beta}f_X^{\textup{MI}}}{\partial x^{\beta}}(x_{\mathcal B}) + h^{\lfloor \alpha \rfloor} \sum_{|\beta|=\lfloor \alpha \rfloor}R_{\beta}(x_{\mathcal B}+{v}h){v}^{\beta}
\end{align}
with 
\begin{align}\label{fMI}
	R_{\beta }(x_{\mathcal B}+{v}h)=\frac {\lfloor \alpha \rfloor }{\beta !}\int _{0}^{1}(1-t)^{\lfloor \alpha \rfloor -1}\frac{\partial^{\beta}f_X^{\textup{MI}}}{\partial x^{\beta}}(x_{\mathcal B}+t {{v}}h)\,dt.
\end{align}
Then, with similar arguments as in the proof of Lemma \ref{lem:norme_sup_bertin}, together with Assumptions \ref{hyp:lipschitzbis} and \ref{hyp:kernelbis} and the fact that $({\mathcal X}_{\mathcal B}')_{\mathcal B \subset \mathcal A}$ forms a partition of $[-1,1)^d$, we get
	\begin{align*}
	|\E&[\widetilde f_{n,h,i}(x)] -f_X(x)| \\& \leqslant 
	3 C_{f_X} \lfloor \alpha \rfloor  h^{\alpha}\, \Bigl(\sum_{|\beta|=\lfloor \alpha \rfloor } \frac{1}{\beta !}\Bigr)  
	\Bigl(\int _{0}^{1}(1-t)^{\lfloor \alpha \rfloor -1} t^{\alpha-\lfloor \alpha \rfloor } dt\Bigr)
	\Bigl(\sum_{\mathcal B \subset \mathcal A}  \int_{ {{\mathcal X}_{\mathcal B}'}}  |{v}|^{\alpha} | {\widetilde K({v})} | d{v}\Bigr) \notag\\
	&\leqslant 3 C_{f_X} \lfloor \alpha \rfloor  h^{\alpha}  \Bigl(\sum_{|\beta|=\lfloor \alpha \rfloor } \frac{1}{\beta !} \Bigr)  \Bigl(\int _{0}^{1}(1-t)^{\lfloor \alpha \rfloor -1} t^{\alpha-\lfloor \alpha \rfloor } dt\Bigr) \Bigl(\int_{[-1,1]^d} | {\widetilde K({v})} | d{v}\Bigr) 
	= C h^{\alpha}. \qedhere
\end{align*}
\end{proof}}

We state Lemma \ref{lem:mise_g_pujol} below as a preliminary result to the proof of  (ii) of Theorem \ref{th:mainresult_compact}.

\begin{lem}\label{lem:mise_g_pujol}
	Let $\alpha >0$. Under Assumptions \ref{hyp:domain}, \ref{hyp:lipschitz_g},  \ref{hyp:kernelbis}, and \ref{hyp:VCclass},
	\begin{align}\label{eq:mise_g_pujol}
		\forall i\in\{1,\cdots,d\},\quad \E\bigl[\int (\widetilde g_{n,h,i}(x)-g(x))^2 dx \bigr]=o(n^{-1/2})
	\end{align}
	as soon as 
	$nh^{2d}\to \infty$, and $nh^{4\alpha}\to 0$. 
\end{lem}

{ 
\begin{proof}[Proof of Lemma \ref{lem:mise_g_pujol}] The proof follows similar lines as the one of Lemma \ref{lem:mise_g_bertin}.
Recall that for $x \in [0,1]^d$, $g(x)=f_X(x)m(x)$.
We extend $m$ as $m^{\textup{MI}}$ on $[-1,2]^d$ as follows. 
 {For all $y \in [-1,2]^d$, 
$m^{\textup{MI}}(y)=m(y^{\textup{MI}})$ with $y^{\textup{MI}}$ defined in the proof of Lemma \ref{lem:norme_sup_pujol}.
Then, for all $x \in [0,1]^d$ and $a \in \{-1,0,1\}^d$, 
$m^{\textup{MI}}(M^a(x))=m(x)$ and similarly, for any $\mathcal B \subset \mathcal A$, $m^{\textup{MI}}(x_{\mathcal B})=m(x)$. We introduce the function $g^{\textup{MI}}(x)=f_X^{\textup{MI}}(x) m^{\textup{MI}}(x)$. We also introduce the intermediate function ${\widetilde g^{\textup{MI}}} (x)=f_X^{\textup{MI}}(x)m(x)$. }
Now following the proof of Lemma \ref{lem:mise_g_bertin}, we get, for all $x \in [0,1]^d$,
$\mathbb{E}[\widetilde{g}_{n,h,i}(x)]-g(x)=\sum_{\mathcal B \subset \mathcal A} \int_{ {{\mathcal X}_{\mathcal B}'}}  { \widetilde K({v})}\Bigl(g^{\textup{MI}}(x_{\mathcal B}+{v}h)-g^{\textup{MI}}(x_{\mathcal B})\Bigr)d{v}$. Now, 
$$\begin{array}{rcl}
	g^{\textup{MI}}(x_{\mathcal B}+{v}h)-g^{\textup{MI}}(x_{\mathcal B})&=&g^{\textup{MI}}(x_{\mathcal B}+{v}h) - {\widetilde g^{\textup{MI}}}(x_{\mathcal B}+{v}h)+{\widetilde g^{\textup{MI}}}(x_{\mathcal B}+{v}h)-{\widetilde g^{\textup{MI}}}(x_{\mathcal B})\\
	& & \quad \quad +{\widetilde g^{\textup{MI}}}(x_{\mathcal B}+{v}h)-{\widetilde g^{\textup{MI}}}(x_{\mathcal B})+{\widetilde g^{\textup{MI}}}(x_{\mathcal B})- g^{\textup{MI}}(x_{\mathcal B})\\
	& =& f_X^{\textup{MI}}(x_{\mathcal B}+{v}h) \Bigl(m^{\textup{MI}}(x_{\mathcal B}+{v}h)-m(x_{\mathcal B}+{v}h)\Bigr)\\
	&& \quad  +{\widetilde g^{\textup{MI}}}(x_{\mathcal B}+{v}h)-{\widetilde g^{\textup{MI}}}(x_{\mathcal B})+f_X^{\textup{MI}}(x_{\mathcal B}) \Bigl(m(x_{\mathcal B})-m^{\textup{MI}}(x_{\mathcal B})\Bigr)\, .
	\end{array}$$ 
Now by definition, we have $m^{\textup{MI}}(x_{\mathcal B}+{v}h)=m(x_{\mathcal B}+{v}h)$ and $m^{\textup{MI}}(x_{\mathcal B})=m(x_{\mathcal B})$. Thus we have $g^{\textup{MI}}(x_{\mathcal B}+{v}h)-g^{\textup{MI}}(x_{\mathcal B})= {\widetilde g^{\textup{MI}}}(x_{\mathcal B}+{v}h)-{\widetilde g^{\textup{MI}}}(x_{\mathcal B})$. As $f_X^{\textup{MI}}$ and $m$ belong to ${\mathcal C}^{\alpha}( [-1,2]^d)$ from Assumptions  \ref{hyp:lipschitzbis} and \ref{hyp:lipschitz_g}, the function ${\widetilde g^{\textup{MI}}}$ also belongs to ${\mathcal C}^{\alpha}( [-1,2]^d)$. 
Then, mimicking \eqref{eq:taylor_fMI} and \eqref{fMI}, we obtain for any $y \in {\mathcal Y }$,
\begin{align*}
{\widetilde g^{\textup{MI}}}(x_{\mathcal B}+{v}h)-{\widetilde g^{\textup{MI}}}(x_{\mathcal B})
	&=\sum_{1\leqslant |\beta|< \lfloor \alpha \rfloor }  \frac{h^{|\beta|}}{\beta \, !} {v}^{\beta} \frac{\partial^{\beta}{\widetilde g^{\textup{MI}}}}{\partial x^{\beta}}(x_{\mathcal B}) + h^{\lfloor \alpha \rfloor} \sum_{|\beta|=\lfloor \alpha \rfloor}R_{\beta}(x_{\mathcal B}+{v}h){v}^{\beta}
\end{align*}
with 
\begin{align*}
	R_{\beta }(x_{\mathcal B}+{v}h)=\frac {\lfloor \alpha \rfloor }{\beta !}\int _{0}^{1}(1-t)^{\lfloor \alpha \rfloor -1}\frac{\partial^{\beta}{\widetilde g^{\textup{MI}}}}{\partial x^{\beta}}(x_{\mathcal B}+t {{v}}h)\,dt.
\end{align*}
Then, we conclude with similar arguments as in the proof of Lemma \ref{lem:norme_sup_bertin}. From Assumption  \ref{hyp:kernelbis} and as $({\mathcal X}_{\mathcal B}')_{\mathcal B \subset \mathcal A}$ forms a partition of $[-1,1)^d$, we have
$$\sum_{\mathcal B \subset \mathcal A} \int_{ {{\mathcal X}_{\mathcal B}'}}  { \widetilde K({v})}\sum_{|\beta|=\lfloor \alpha \rfloor}{v}^{\beta} d{v} =\int_{-1}^1 \Bigl(\sum_{|\beta|=\lfloor \alpha \rfloor}{v}^{\beta}\Bigr)\widetilde K({v}) d{v} =0 \, .$$
Thus 
$$\begin{array}{rcl}
	\mathbb{E}[\widetilde{g}_{n,h,i}(x)]-g(x)& = & \sum_{\mathcal B \subset \mathcal A} \int_{ {{\mathcal X}_{\mathcal B}'}}  { \widetilde K({v})}\Bigl(g^{\textup{MI}}(x_{\mathcal B}+{v}h)-g^{\textup{MI}}(x_{\mathcal B})\Bigr)d{v}\\
	& = & \sum_{\mathcal B \subset \mathcal A} \int_{ {{\mathcal X}_{\mathcal B}'}}  { \widetilde K({v})}h^{\lfloor \alpha \rfloor} \sum_{|\beta|=\lfloor \alpha \rfloor}R_{\beta}(x_{\mathcal B}+{v}h){v}^{\beta}d{v} \\
	& = &  \sum_{\mathcal B \subset \mathcal A} \int_{ {{\mathcal X}_{\mathcal B}'}}  { \widetilde K({v})}h^{\lfloor \alpha \rfloor} \sum_{|\beta|=\lfloor \alpha \rfloor}\frac {\lfloor \alpha \rfloor }{\beta !}\\
	& & \quad \quad \int _{0}^{1}(1-t)^{\lfloor \alpha \rfloor -1}\Bigl(\frac{\partial^{\beta}{\widetilde g^{\textup{MI}}}}{\partial x^{\beta}}(x_{\mathcal B}+t {{v}}h)-\frac{\partial^\beta {\widetilde g^{\textup{MI}}}}{\partial x^\beta}(x_{\mathcal B})\Bigr)dt {v}^{\beta}d{v}\, .
	\end{array}$$
Now, using \eqref{eq:holder} since  ${\widetilde g^{\textup{MI}}}\in\mathcal C^{\alpha}([0,1]^d)$, we have 
\begin{align*}
	\Bigl |\frac{\partial^\beta {\widetilde g^{\textup{MI}}}}{\partial x^\beta}(x_{\mathcal B}+t {{v}}h) - \frac{\partial^\beta {\widetilde g^{\textup{MI}}}}{\partial x^\beta}(x_{\mathcal B})\Bigr |
	\leqslant C( h t\|{v} \|_{\infty})^{\alpha - \lfloor \alpha \rfloor}\leqslant C h^{\alpha - \lfloor \alpha \rfloor}
\end{align*}
for all ${v} \in [0,1]^d$, $x\in [0,1]^d$, $t\in [0,1]$, $h\in (0,\infty)$, and $\beta \in \N^d$ such that $|\beta| = \lfloor \alpha\rfloor$. We hence conclude that 
$\displaystyle \mathbb{E}[\widetilde{g}_{n,h,i}(x)]-g(x)\leq C h^{\alpha} $.
\end{proof}
}

\subsection{Proof of Theorem \ref{th:mainresult_compact}}\label{app:mainresult_compact}

\begin{proof}[Proof of Theorem \ref{th:mainresult_compact}]
Following the same lines as in the proof of Theorem 2.1 in \cite{doksum1995nonparametric}, we aim at proving that 
\begin{align}\label{eq:a_prouver}
\widehat T_{n,h}-\E[\E[Y|X]^2]=\frac 1n \sum_{i=1}^n (2Y_i-m_i)m_i 
-\E[\E[Y|X]^2]+o_\P(n^{-1/2}).
\end{align}
The conclusion of Theorem \ref{th:mainresult_compact} will then follow directly applying the standard central limit theorem for the sum of i.i.d.\ random variables to the right-hand side of the previous display together with Slutsky's lemma. To establish \eqref{eq:a_prouver}, we compute
\begin{align*}
\widehat T_{n,h}-\frac 1n \sum_{i=1}^n (2Y_i-m_i)m_i 
&= \frac 1n \sum_{i=1}^n \bigl[(2Y_i -\widehat m_i)\widehat m_i - (2Y_i - m_i) m_i\bigr]\\
&= \frac 1n \sum_{i=1}^n \bigl[2Y_i (\widehat m_i-m_i)+ m_i^2-\widehat m_i^2\bigr]\\
&= \frac 1n \sum_{i=1}^n \bigl[2(Y_i-m_i) (\widehat m_i-m_i)-( \widehat m_i- m_i)^2\bigr]\eqdef I_1-I_2.
\end{align*}

\paragraph{Study of $I_1$}
Since $\widehat m_i-m_i = \frac{\widehat g_i-\widehat f_im_i}{f_i} + \frac{(f_i-\widehat f_i)(\widehat g_i-\widehat f_im_i)}{f_i\widehat f_i}$, $I_1$ rewrites as the sum of two terms $I_{11}$ and $I_{12}$ with 
\begin{align*}
I_{11}& = \frac 2n \sum_{i=1}^n (Y_i-m_i) \frac{\widehat g_i-\widehat f_im_i}{f_i}\\
&= \frac {2}{n(n-1)} \sum_{i=1}^n \sum_{\substack{j=1,\\ j\neq i}}^n (Y_i-m_i)(Y_j-m_i) \frac{K_h\circ A_{X_i}(X_j-X_i)}{f_i}\\
&= \frac {2}{n(n-1)} \sum_{i=1}^n \sum_{\substack{j=1,\\ j\neq i}}^n \varepsilon_i(\varepsilon_j+m_j-m_i) \frac{K_h\circ A_{X_i}(U_{ij})}{f_i}
\end{align*}
denoting the residual $Y_i-m_i$ by $\varepsilon_i$ and the difference $X_j-X_i$ by $U_{ij}$. Conditioning by $X^{(n)}=(X_1,\cdots,X_n)$ and using independence, we deduce that
\begin{align*}
\E[I_{11}]&=2\E\Bigl[\varepsilon_1(\varepsilon_2+m_2-m_1) \frac{K_h\circ A_{X_1}(U_{12})}{f_i}\Bigr]\\
&= 2\E\Bigl[\E[\varepsilon_1\vert X_1] (\E[\varepsilon_2\vert X_2]+m_2-m_1) \frac{K_h\circ A_{X_1}(U_{12})}{f_i}\Bigr]
\end{align*}
that cancels since $E[\varepsilon_1\vert X_1]=0$ while $\E[I_{11}^2]$ equals
\begin{align*} 
&\frac{4}{n^2(n-1)^2} \sum_{i,k=1}^n \sum_{\substack{j,\ell=1,\\ j\neq i, \\ \ell\neq k}}^n 
\E\Bigl[\varepsilon_i \varepsilon_k (\varepsilon_j+m_j-m_i) (\varepsilon_\ell+m_\ell-m_k) \frac{K_h\circ A_{X_i}(U_{ij})}{f_i} \frac{K_h\circ A_{X_k}(U_{k\ell})}{f_k} \Bigr]\\
&= \frac{4}{n^2(n-1)^2} \sum_{i=1}^n \sum_{\substack{j=1,\\ j\neq i}} 
\E\Bigl[\varepsilon_i^2(\varepsilon_j+m_j-m_i)^2  \frac{(K_h\circ A_{X_i}(U_{ij}))^2}{f_i^2} \Bigr]\\
&\quad + \frac{4}{n^2(n-1)^2} \sum_{i=1}^n \sum_{\substack{j=1,\\ j\neq i}} 
\E\Bigl[\varepsilon_i \varepsilon_j (\varepsilon_j+m_j-m_i) (\varepsilon_i+m_i-m_j) \frac{K_h\circ A_{X_i}(U_{ij})}{f_i} \frac{K_h\circ A_{X_j}(U_{ji})}{f_j} \Bigr].
\end{align*}
Conditioning by $X^{(n)}=(X_1,\cdots,X_n)$ once again leads to
\begin{align*}
\E[\varepsilon_i^2&(\varepsilon_j+m_j-m_i)^2\vert X^{(n)}]
=\E[\varepsilon_i^2\varepsilon_j^2\vert X^{(n)}]+ 2\E[\varepsilon_i^2\varepsilon_j(m_j-m_i)\vert X^{(n)}]+ \E[\varepsilon_i^2(m_j-m_i)^2\vert X^{(n)}]\\
&=\E[\varepsilon_i^2\vert X_i]\E[\varepsilon_j^2\vert X_j]+ 2(m_j-m_i) \E[\varepsilon_i^2\vert X_i]\E[\varepsilon_j\vert X_j]+ (m_j-m_i)^2\E[\varepsilon_i^2\vert X_i]\\
&=\sigma_i^2\sigma_j^2+(m_j-m_i)^2 \sigma_i^2
\end{align*}
and
\begin{align*}
&\E[\varepsilon_i\varepsilon_j(\varepsilon_j+m_j-m_i)(\varepsilon_i+m_i-m_j)\vert X^{(n)}]\\
&=\E[\varepsilon_i^2\varepsilon_j^2\vert X^{(n)}]+ (m_j-m_i)\E[\varepsilon_i^2\varepsilon_j\vert X^{(n)}]+
(m_i-m_j)\E[\varepsilon_j^2\varepsilon_i\vert X^{(n)}] -(m_j-m_i)^2\E[\varepsilon_i\varepsilon_j\vert X^{(n)}]\\
&=\sigma_i^2\sigma_j^2.
\end{align*}
Hence
\begin{align*}
\E[I_{11}^2] &= \frac{4}{n^2(n-1)^2} \sum_{i=1}^n \sum_{\substack{j=1,\\ j\neq i}} 
\E\Bigl[(\sigma_i^2\sigma_j^2+(m_j-m_i)^2 \sigma_i^2)  \frac{(K_h\circ A_{X_i}(U_{ij}))^2}{f_i^2} \Bigr]\\
&\qquad + \frac{4}{n^2(n-1)^2} \sum_{i=1}^n \sum_{\substack{j=1,\\ j\neq i}} 
\E\Bigl[\sigma_i^2\sigma_j^2 \frac{K_h\circ A_{X_i}(U_{ij})}{f_i} \frac{K_h\circ A_{X_j}(U_{ji})}{f_j} \Bigr].
\end{align*}

By Assumptions \ref{hyp:densityb}, \ref{hyp:bounds}, \ref{hyp:kernel}, and the fact that $m$ is bounded from above, we conclude that
\begin{align*}
\E[I_{11}^2] &= O\Bigl(\frac{1}{n^2h^{2d}}\Bigr)=o\Bigl(\frac{1}{n}\Bigr)
\end{align*}
since $nh^{2d}\to \infty$ by assumption.

Let us turn to the computation of $I_{12}$ defined as follows 
\begin{align*}
I_{12}& = \frac 2n \sum_{i=1}^n (Y_i-m_i) \frac{(f_i-\widehat f_i)(\widehat g_i-\widehat f_im_i)}{f_i\widehat f_i}.
\end{align*}
By \eqref{eq:hat_f_bounded_bertin} and since $\widehat g_1$ does not depend on { $Y_1$},
\begin{align*}
&\abs{\E[I_{12}]}\leqslant \frac {C}{n} \sum_{i=1}^n \E_{\substack{X_1,\cdots,X_n\\Y_1,\cdots,Y_n}}\Bigl[|\varepsilon_i|\frac{|f_i-\widehat f_i|(|g_i-\widehat g_i|+|f_i-\widehat f_i|| m_i|)}{f_i}\Bigr]\\
&=C~\E_{\substack{X_1,\cdots,X_n\\Y_1,\cdots,Y_n}}\Bigl[|\varepsilon_1|\frac{|f_1-\widehat f_1|(|g_1-\widehat g_1|+|f_1-\widehat f_1 ||m_1|)}{f_1}\Bigr]\\
&=C~\E_{X_1,\cdots,X_n}\Bigl[\frac{|f_1-\widehat f_1|}{f_1} 
\E_{Y_1,\cdots,Y_n}[|\varepsilon_1|(|g_1-\widehat g_1|+|f_1-\widehat f_1|| m_1|)\vert X^{(n)}]\Bigr]\\
&=C~\E_{X_1,\cdots,X_n}\Bigl[\frac{|f_1-\widehat f_1|}{f_1} 
\E_{Y_1}[|\varepsilon_1|\vert X_1](\E_{Y_2,\cdots,Y_n}[|g_1-\widehat g_1|\vert X^{(n)}]+|f_1-\widehat f_1|| m_1|)\Bigr]\\
&\leqslant C~\E_{X_1,\cdots,X_n}\Bigl[\frac{|f_1-\widehat f_1|}{f_1} 
\sigma_1(\E_{Y_2,\cdots,Y_n}[|g_1-\widehat g_1|\vert X^{(n)}]+|f_1-\widehat f_1|| m_1|)\Bigr]\\
&= C~\E_{\substack{X_1,\cdots,X_n\\Y_2,\cdots,Y_n}} \Bigl[\frac{|f_1-\widehat f_1|}{f_1} 
\sigma_1(|g_1-\widehat g_1|+|f_1-\widehat f_1|| m_1|)\Bigr]\\
&= C~\E_{\substack{X_2,\cdots,X_n\\Y_2,\cdots,Y_n}}\Bigl[\int_{[0,1]^d}  |f_X(x)-\widehat f_{n,h,1}(x)| 
\sigma(x)(|g(x)-\widehat g_{n,h,1}(x)|+|f_X(x)-\widehat f_{n,h,1}(x)|| m(x)|)dx\Bigr]\\
&\leqslant C~\E_{\substack{X_2,\cdots,X_n\\Y_2,\cdots,Y_n}}\Bigl[\int_{[0,1]^d}  |f_X(x)-\widehat f_{n,h,1}(x)| 
|g(x)-\widehat g_{n,h,1}(x)|dx\Bigr]\\
&\qquad +C~\E_{\substack{X_2,\cdots,X_n\\Y_2,\cdots,Y_n}}\Bigl[\int_{[0,1]^d}  |f_X(x)-\widehat f_{n,h,1}(x)||f_X(x)-\widehat f_{n,h,1}(x)|| m(x)|)dx\Bigr].
\end{align*}
By applying twice Cauchy-Schwartz inequality, the first term is bounded from above by
\begin{align*}
C &~\E_{\substack{X_2,\cdots,X_n\\Y_2,\cdots,Y_n}}\Bigl[\bigl(\int_{[0,1]^d}  ( f_X(x)-\widehat f_{n,h,1}(x))^2 dx \bigr)^{1/2} 
\bigl(\int_{[0,1]^d}  (g(x)-\widehat g_{n,h,1}(x))^2dx\bigr)^{1/2}\Bigr]\\
&\leqslant C~\E_{\substack{X_2,\cdots,X_n\\Y_2,\cdots,Y_n}}\Bigl[\int_{[0,1]^d}  ( f_X(x)-\widehat f_{n,h,1}(x))^2 dx\Bigr] ^{1/2}
\E_{\substack{X_2,\cdots,X_n\\Y_2,\cdots,Y_n}}\Bigl[\int_{[0,1]^d}  (g(x)-\widehat g_{n,h,1}(x))^2dx\Bigr]^{1/2}\\
&=o(n^{-1/2})
\end{align*}
using \eqref{eq:mise_f_bertin} and \eqref{eq:mise_g_bertin} while the second term is $o(n^{-1/2})$.

\paragraph{Study of $I_2$} Recall that
\begin{align*}
I_2&= \frac 1n \sum_{i=1}^n ( \widehat m_i- m_i)^2=\frac 1n \sum_{i=1}^n \frac{(\widehat g_i- \widehat f_i m_i)^2}{\widehat f_i^2}= { O_{\mathbb{P}}\Bigl(\frac{1}{n} \sum_{i=1}^n (\widehat g_i- \widehat f_i m_i)^2\Bigr)}
\end{align*}
by \eqref{eq:hat_f_bounded_bertin}. Then,
\begin{align*}
&\E[I_2] \leqslant C~\E_{\substack{X_2,\cdots,X_n\\Y_2,\cdots,Y_n}}\Bigl[\int_{[0,1]^d}  (\widehat g_{n,h,i}(x)- \widehat f_{n,h,i}(x) m(x))^2 f_X(x)dx\Bigr]\\
&= C~\E_{\substack{X_2,\cdots,X_n\\Y_2,\cdots,Y_n}}\Bigl[\int_{[0,1]^d}  [(\widehat g_{n,h,i}(x)- g(x))+(f_X(x)m(x)-\widehat f_{n,h,i}(x) m(x))]^2 f_X(x)dx\Bigr]\\
&\leqslant C~\E_{\substack{X_2,\cdots,X_n\\Y_2,\cdots,Y_n}}\Bigl[\int_{[0,1]^d}  (\widehat g_{n,h,i}(x)-g(x))^2 f_X(x)dx\Bigr]+C~\E\Bigl[\int_{[0,1]^d} ( \widehat f_{n,h,i}(x)-f(x))^2 m(x)^2f_X(x)dx\Bigr]\\
&=o(n^{-1/2})
\end{align*}
as $f_X$ is continuous thus bounded on $[0,1]^d$, as $m$ is also bounded on $[0,1]^d$ (recall that it follows from the boundedness of $\sigma^2$) and from the results stated in \eqref{eq:mise_f_bertin} and  \eqref{eq:mise_g_bertin}.

\medskip

As for item (ii), the proof is similar  except that for any $i \neq j$, $K_h\circ A_{X_i}(U_{ij})$ is replaced by $\sum_{a \in \{-1,0,1\}^d}\prod_{k=1}^d K_h(M^{a_k}(X_{k,j})-X_i)$.
\end{proof}

\subsection{Proof of the remaining results}\label{app:maincor_compact}

\begin{proof}[Proof of Proposition \ref{prop:asymp_eff}]
We denote by $P$ the distribution of $(X,Y)$ and we introduce 
\[
\psi(P)=\E[\E[Y\vert X]^2].
\]
The influence efficient function of $\psi$ at $P$, as stated in \cite {doksum1995nonparametric}, is given by  
$\widetilde \psi_P(x,y)=(2y-m(x))m(x) - \E[Ym(X)]$ (see \cite{klein2024note_efficiency} for explicit computations). 
Moreover, we deduce from \eqref{eq:a_prouver} that
\[
\widehat T_{n,h} = \psi(P)+ \frac 1 n \sum_{i=1}^n \widetilde \psi_{P}(X_i,Y_i) + o_\P(n^{-1/2})
\]
and, by \cite[Condition (25.22)]{van2000asymptotic}, $\widehat T_{n,h}$ is asymptotically efficient and so is $\widetilde T_{n,h}$. 
\end{proof}

\begin{proof}[Proof of Corollary \ref{cor:TCL_global}]

In view of the definition of the asymptotic efficiency (see, e.g., \cite[Lemma 25.23]{van2000asymptotic} or \cite[Definition 2.4]{klein2024note_efficiency}, the asymptotic efficiency ensures the asymptotic normality. 
Thus we only need to prove the asymptotic efficiency. To do so, it suffices to use the asymptotic efficiency of $\widehat T_{n,h}$ and $\widetilde T_{n,h}$ in Proposition \ref{prop:asymp_eff}, the asymptotic efficiency of the empirical mean $\overline Y_n$ and the empirical variance $\overline{Y^2}_n$ (see, e.g., \cite{janon2012asymptotic}) together with the efficiency in product space \cite[Theorem 25.50]{van2000asymptotic} to conclude to the joint asymptotic efficiency of 
$(\widehat T_{n,h}, \overline Y_n,\overline Y_n^2)$ and $(\widetilde T_{n,h}, \overline Y_n,\overline Y_n^2)$. Finally, we follow the same lines as in the proof of \cite[Proposition 2.5]{janon2012asymptotic}
using the efficiency and delta method \cite[Theorem 25.47]{van2000asymptotic} to get the required result.  It remains to perform easy computations to get the expression of the asymptotic variance
{
\begin{align}
\Var  (Y)^2  \sigma^2
 = &~~ \Var((2Y-m(X))m(X)) + 
   S^X (\Var(Y^2)S^X
  -2 \Cov(m(X)^2,Y^2) )\notag\\
&   +  4 \E[Y] (S^X-1) [\E[m(X)^2Y]- \E[Y] \E[Y^2] - \Cov(Y,Y^2)
S^X ],\label{eq:sigma2}
\end{align}
from which we derive simpler expressions 
\begin{itemize}
    \item when $\E[Y]=0$ or $S^X=1$: 
\begin{equation*}
\Var(Y)^2 \sigma^2
 =  \Var((2Y-m(X))m(X)) + S^X (\Var(Y^2)S^X
  -2 \Cov(m(X)^2,Y^2) ),
 \end{equation*}
\item when $S^X=0$: 
\begin{equation*}
\Var(Y)^2 \sigma^2
 =  \Var((2Y-m(X))m(X))
  + 4 \E[Y] (S^X-1) [\E[m(X)^2Y]- \E[Y] \E[Y^2]].
\end{equation*}
\end{itemize}
}
The proof is then complete. 
\end{proof}


\newpage

\section{Additional numerical experiments}\label{app:xps}

\subsection{{Influence of $\varepsilon$ in \cite{doksum1995nonparametric}}}

\begin{figure}[h]
\centering
    \includegraphics[width=\textwidth]{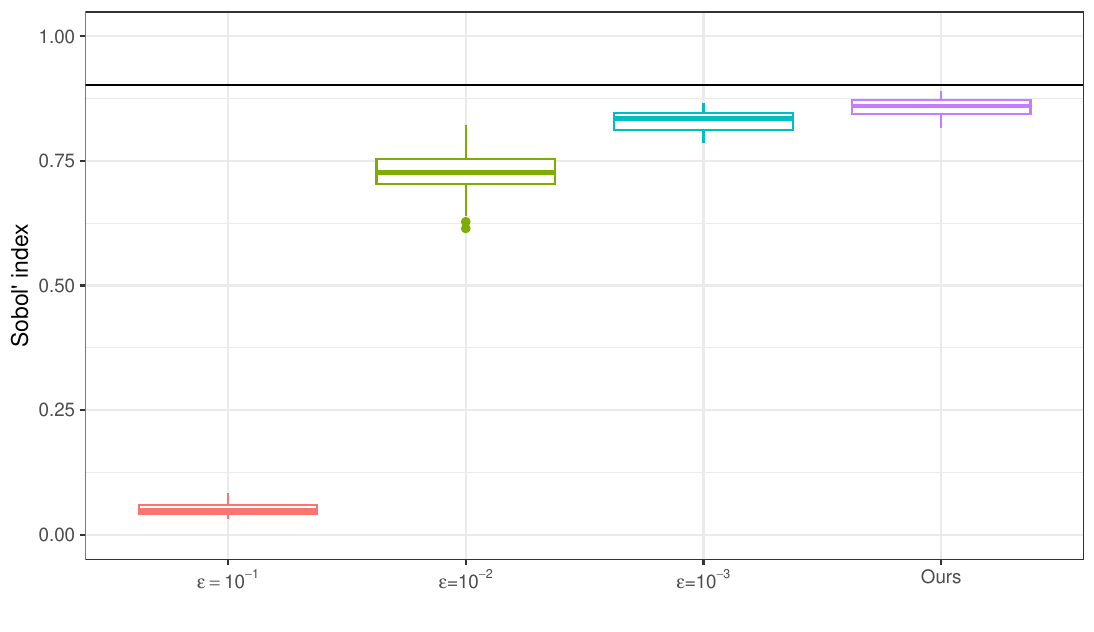}
\caption{Comparison of our mirror-type estimator with the estimator of \cite{doksum1995nonparametric} for different values of $\varepsilon$. The reference value of the index is represented with a gray line.}
\label{fig:DK}
\end{figure}

We illustrate numerically that the choice of the $\varepsilon$ tuning parameter of the estimator proposed in \cite{doksum1995nonparametric} is very sensitive, thus limiting its practical use as opposed to our mirror-type estimator. We consider Example 3.2 from \cite{doksum1995nonparametric} and test $\varepsilon=10^{-1},10^{-2},10^{-3}$. The comparison with our estimator with a kernel of order 2 is given in Figure \ref{fig:DK}. When $\varepsilon$ is equal to $10^{-3}$, the performance of both estimators are similar. However in other cases, the bias of  \cite{doksum1995nonparametric} can be very large. Since in practice such an estimation problem is unsupervised, the tuning of $\varepsilon$ seems highly difficult and the non-robustness of the final estimator with respect to this parameter limits its practical use.

\subsection{Analytical test functions with smaller sample size}

For completeness, we include below new results on the Bratley and g-Sobol functions in dimension $p=5$ with smaller sample sizes. As expected, the bias is larger when $n$ decreases.

\begin{figure}[H]
\centering
\begin{subfigure}[b]{0.45\textwidth}
    \centering
    \includegraphics[width=\textwidth]{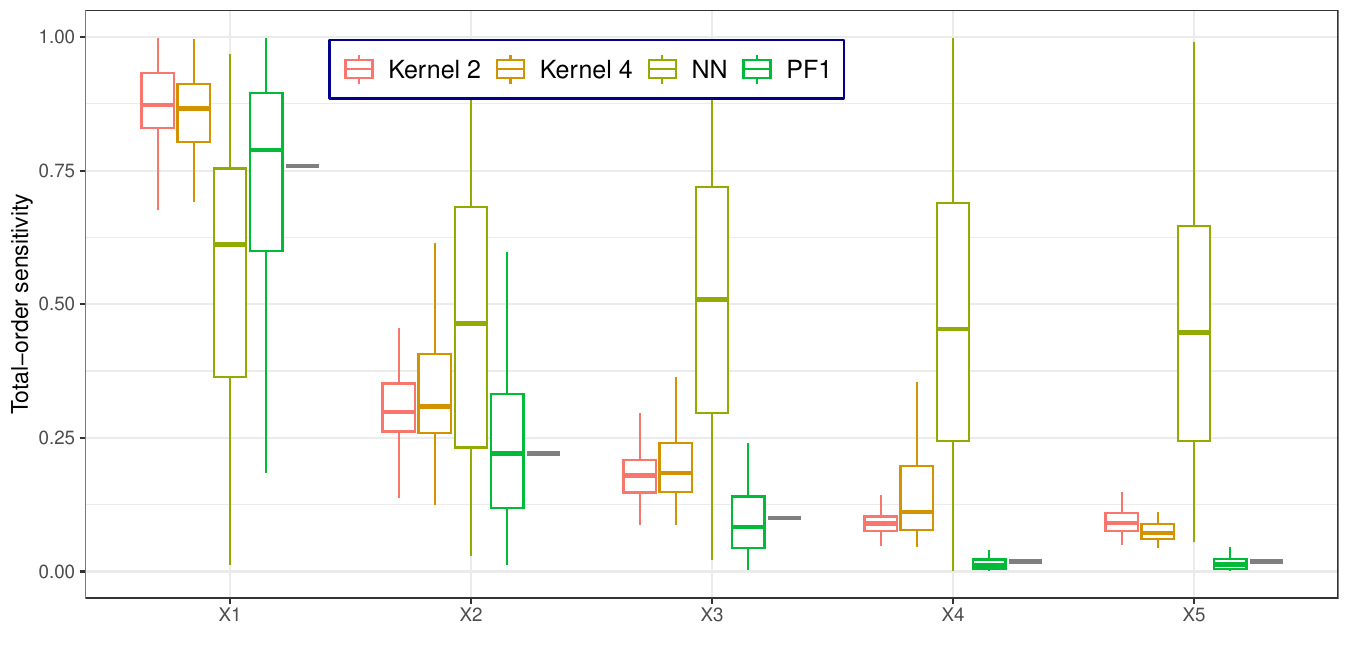}
    \caption{$n=100$}
\end{subfigure}
\hfill
\begin{subfigure}[b]{0.45\textwidth}
    \centering
    \includegraphics[width=\textwidth]{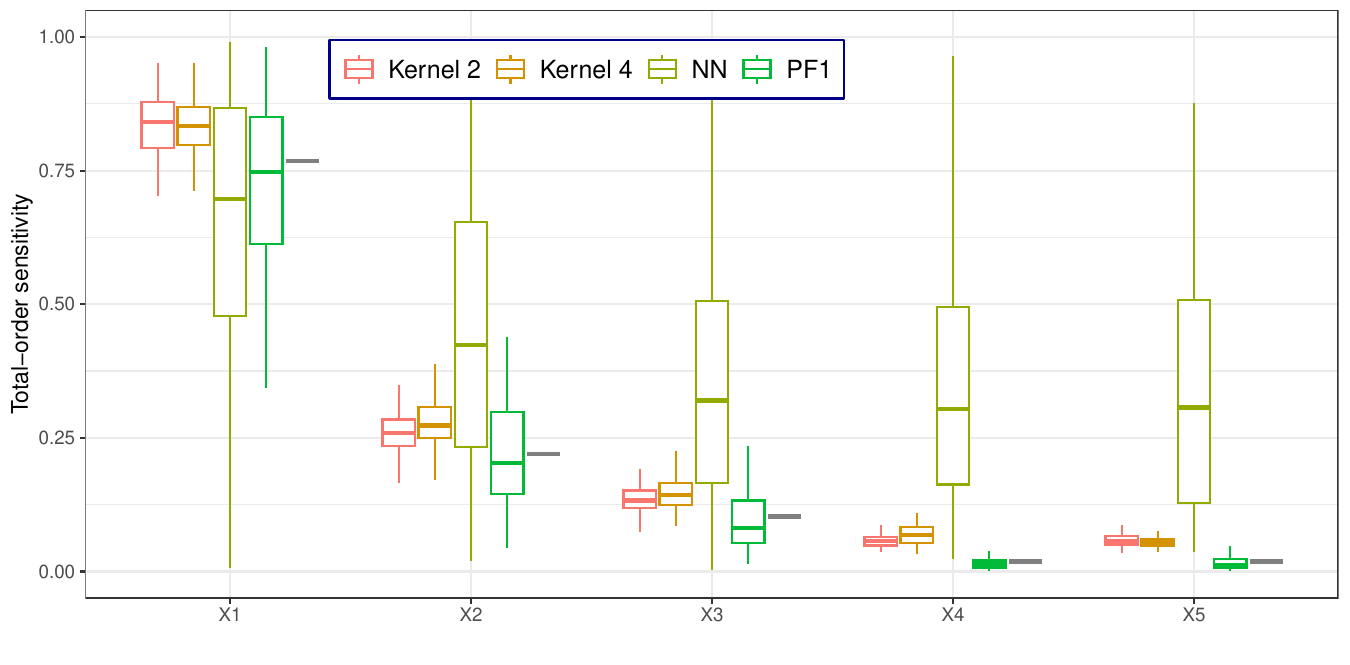}
    \caption{$n=200$}
\end{subfigure}
\begin{subfigure}[b]{0.45\textwidth}
    \centering
    \includegraphics[width=\textwidth]{Bratley_dim_5_n_500_rep_100_eta1_Boxplot_T_final.pdf}
    \caption{$n=500$}
\end{subfigure}
\hfill
\begin{subfigure}[b]{0.45\textwidth}
    \centering
    \includegraphics[width=\textwidth]{Bratley_dim_5_n_1000_rep_100_eta1_Boxplot_T_final.pdf}
    \caption{$n=1000$}
\end{subfigure}
\caption{Estimators for total indices of the Bratley function from $n=100$ to $n=1000$. The reference value is represented with a gray line.}
\end{figure}

\begin{figure}[H]
\centering
\begin{subfigure}[b]{0.45\textwidth}
    \centering
    \includegraphics[width=\textwidth]{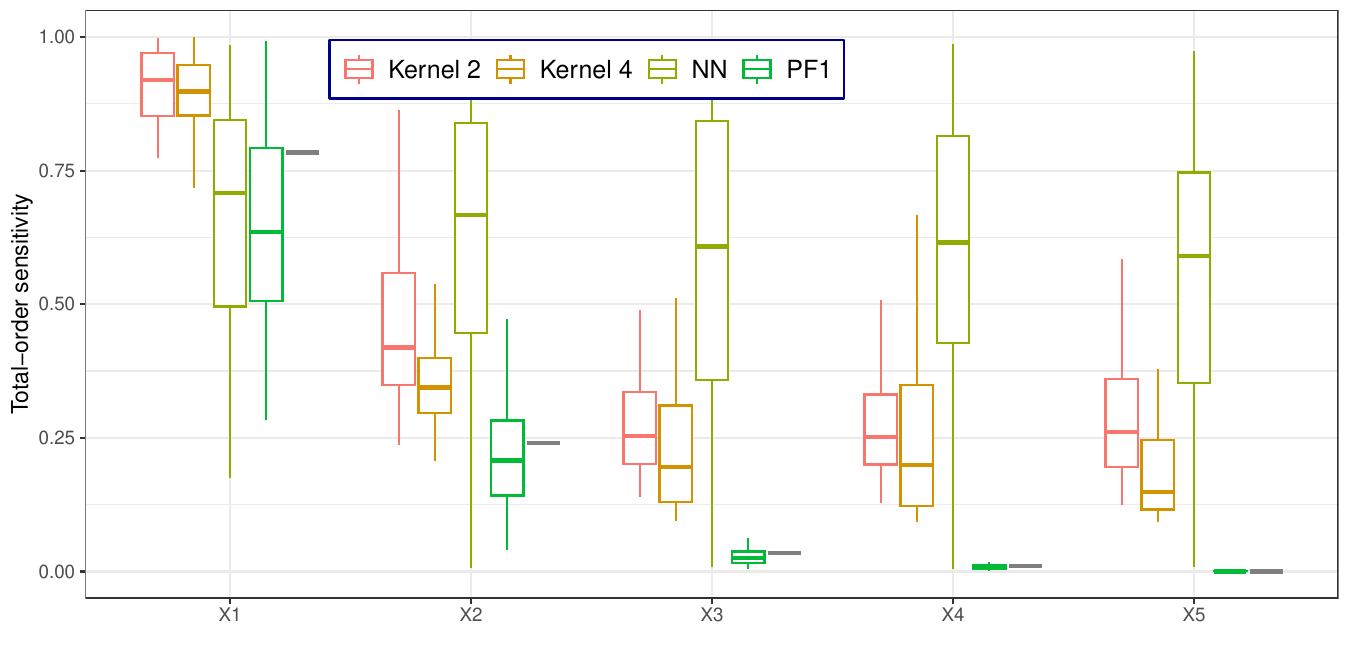}
    \caption{$n=100$}
\end{subfigure}
\hfill
\begin{subfigure}[b]{0.45\textwidth}
    \centering
    \includegraphics[width=\textwidth]{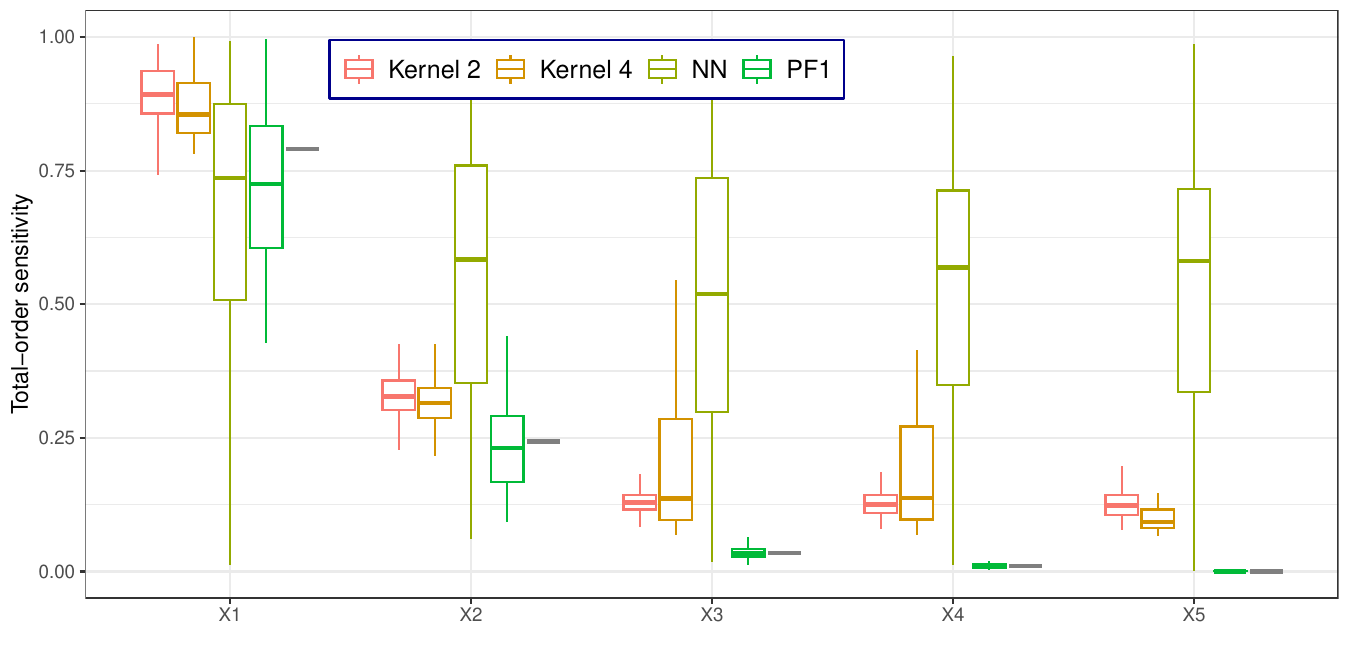}
    \caption{$n=200$}
\end{subfigure}
\begin{subfigure}[b]{0.45\textwidth}
    \centering
    \includegraphics[width=\textwidth]{GSobol_dim_5_n_500_rep_100_eta1_Boxplot_T_final.pdf}
    \caption{$n=500$}
\end{subfigure}
\hfill
\begin{subfigure}[b]{0.45\textwidth}
    \centering
    \includegraphics[width=\textwidth]{GSobol_dim_5_n_1000_rep_100_eta1_Boxplot_T_final.pdf}
    \caption{$n=1000$}
\end{subfigure}
\caption{Estimators for total indices of the g-Sobol function from $n=100$ to $n=1000$. The reference value is represented with a gray line.}
\end{figure}

\newpage

\textbf{Acknowledgement} 

\medskip

The authors would like to thank Elmar Plischke for pointing out some very relevant references related to the present work. They are also grateful to the reviewers for their valuable comments that allowed to improve this article.  Support from the ANR-3IA Artificial and Natural Intelligence Toulouse Institute, the ANR GATSBII (ANR-24-CE23-6645) and the consortium in Applied Mathematics CIROQUO-2,
gathering partners in technological research and academia in the development of advanced methods for
Computer Experiments, are gratefully acknowledged.

\bibliographystyle{abbrv}
\bibliography{biblio_SA2}

\end{document}